\documentclass[11pt,reqno]{amsart}
\usepackage{t1enc}
\usepackage[latin1]{inputenc}
\usepackage{amsmath,latexsym,amssymb,graphicx,dsfont,amsthm,amsfonts}
\usepackage[ps2pdf]{hyperref}
\usepackage{color}
\usepackage{umoline}
\usepackage{verbatim}

\newtheorem{Th}{Theorem}[section]
\newtheorem{Lemma}[Th]{Lemma}
\newtheorem{Cor}[Th]{Corollary}
\newtheorem{Prop}[Th]{Proposition}

\newcommand{\Prob}{\mathbb{P}}
\renewcommand{\P}{\mathbb{P}}
\newcommand{\calL}{\mathcal{L}}

\newcommand{\calI}{\mathcal{I}}

\newcommand{\X}{\mathcal{X}}

\newcommand{\T}{\mathcal{T}}
\newcommand{\N}{\mathbb{N}}
\newcommand{\Z}{\mathbb{Z}}
\newcommand{\ro}{\textbf{r}}
\newcommand{\eps}{\varepsilon}

\renewcommand{\k}{\mathbf{k}}

\newtheorem{Rem}[Th]{Remark}
\newenvironment{Example}{\textbf{Example \arabic{section}.\arabic{Th}:}}

\numberwithin{equation}{section}

\begin{document}

\title{Branching Random Walks on Free Products of Groups}

\author{Elisabetta Candellero, Lorenz A. Gilch and Sebastian M\"uller}

\address{Elisabetta Candellero: Institut f\"ur mathematische Strukturtheorie (Math. C), Graz University of Technology, Steyrergasse 30, A-8010 Graz, Austria}
\address{Lorenz Gilch: Universit\'e de Gen\`eve, Section de Math\'ematiques, 2-4, rue du Li\`evre, 1211 Gen\`eve 4, Switzerland}
\address{Sebastian M\"uller: LATP, Université de Provence, 39, rue F. Joliot Curie,
13453 Marseille, France }

\email{candellero@TUGraz.at, gilch@TUGraz.at, mueller@cmi.univ-mrs.fr}
\urladdr{http://www.math.tugraz.at/$\sim$candellero/, http://www.math.tugraz.at/$\sim$gilch/, \newline http://www.latp.univ-mrs.fr/$\sim$mueller}
\date{\today}

\maketitle



\begin{abstract}

We study certain phase transitions of branching random walks (BRW) on Cayley
graphs of free products. The aim of this paper is to compare the size and
structural properties of the trace, i.e., the subgraph that consists of all
edges and vertices that were visited by some particle, with those of the
original Cayley graph. We investigate the phase when the growth parameter
$\lambda$ is small enough such that the process survives but the trace is not
the original graph. A first result is that the box-counting dimension
of the boundary of the trace exists, is almost surely constant and  equals the Hausdorff dimension which we denote by $\Phi(\lambda)$. The
main result states that the function $\Phi(\lambda)$ has only one point of
discontinuity which is at $\lambda_{c}=R$ where $R$ is the radius of
convergence of the Green function of the underlying random walk. Furthermore,
$\Phi(R)$ is bounded by one half the Hausdorff dimension of the boundary of the original
Cayley graph and the behaviour of $\Phi(R)-\Phi(\lambda)$ as $\lambda \uparrow R$ is classified.

In the case of free products of infinite groups the
end-boundary can be decomposed into words of finite and words of infinite
length. We prove the existence of a phase transition such that if $\lambda\leq
\tilde\lambda_{c}$ the end boundary of the trace consists only of infinite
words and if $\lambda>\tilde\lambda_{c}$ it also contains finite words. In the
last case,  the Hausdorff dimension of the set of ends (of the
trace and the original graph) induced by finite words is strictly smaller than
the one of the ends induced by infinite words.

\end{abstract}

\section{Introduction}

A branching random walk (BRW) is a \emph{growing cloud} of particles that move on an underlying graph $\mathcal{X}$ in discrete time. The process starts with one particle in the root $e$ of the graph. Then at each discrete time step a particle produces offspring particles according to some offspring distribution with mean $\lambda>1$, and then each descendent moves one step according to a random walk on $\mathcal{X}$. Particles branch and move independently of the other particles and the history of the process. 
A first natural question is to ask whether the process eventually \emph{fills up} the whole graph, that is, if every finite subset will eventually be occupied or free of particles. If the BRW visits the whole graph it is called recurrent and transient otherwise. 
As a consequence of Kesten's amenability criterion any BRW  is
recurrent on the Cayley graph of an amenable group. Furthermore, one observes a phase transition on non-amenable groups;  there exists some $\lambda_c>1$ such that a BRW
with $\lambda\leq \lambda_c$ is transient, while it is recurrent otherwise. In the transient case the trace of the BRW, that is, the subgraph that
consists of all edges and vertices that were visited by the BRW, is a proper
random subgraph of the original Cayley graph.  Benjamini and M\"uller \cite{benjamini:10} studied first general qualitative statements of the trace of BRW on groups. In particular, they proved exponential volume growth of the trace in general. However, their approach is rather qualitative and gives no quantitative results on the growth rate. In this article we
study BRW on free products of groups and obtain a precise formula for the growth rate and dimensions of the end boundary of the trace.
One motivation to study BRW on this class of structures lies mainly in the fact that they are among the simplest non-amenable groups. This makes them to a reference and starting point for more complicated non-amenable structures such as, for instance, groups with infinitely many ends or  hyperbolic groups. Besides this, free products of groups are interesting   on their own since they play an important role in some fields of algebraic topology and  in Basse--Serre theory.

The starting point of the present investigation of branching random walks was
the work of Hueter and Lalley \cite{lalley2000}, who studied BRW on homogeneous
trees. We remark that in their setting and notation weak survival is equivalent to transience in our language. In the transient regime the BRW eventually vacates every finite subset and the particle trails converge to the geometric end boundary $\Omega$ of the tree. The \emph{limit set} $\Lambda$ of the BRW is the random subset of the boundary that consists of all ends, where the BRW accumulates. By this we mean  that  each neighbourhood of an end in $\Lambda$ is visited infinitely often by the process. Equivalently, we can define $\Lambda$ as the geometric end boundary of the trace. 

Typical ways of measuring the size of boundaries are  by use of the box-counting dimension (also known as Minkowski dimension) or the Hausdorff dimension. 
In \cite{lalley2000} a formula for the Hausdorff dimension of $\Lambda$ is given for BRW on homogeneous trees. In particular,  it is shown there that the limit set has Hausdorff dimension no larger than one half the Hausdorff dimension of the entire boundary $\Omega$. We extend these results to BRW on free products of groups. 
We prove existence of the box-counting dimension, show that the Hausdorff
dimension equals the box-counting dimension and  present a formula in terms
of generating functions of the underlying random walk, see Theorem \ref{thm:boxdim}.
In the same way we obtain a formula for the Hausdorff dimension of the whole
space of ends, see Theorem \ref{thm:boxdim-omega}. This eventually leads to the
result that the Hausdorff dimension of $\Lambda$ is not larger than one half the Hausdorff dimension of the entire boundary.  Another consequence of the formula of the Hausdorff dimension is that  the dimension varies continuously in the subcritical regime, see Theorem \ref{Th:lambda-function}. This affirms the conjecture made in \cite{benjamini:10} for general non-amenable groups that 
the Hausdorff dimension of the limit set is continuous for $\lambda \neq \lambda_{c}$ and discontinuous at $\lambda_{c}$.  As pointed out in \cite{lalley2000} the very same phenomenon hold for other growth processes (e.g. hyperbolic  branching Brownian motion, isotropic contact process on homogeneous trees)  that exhibit a phase transition between \emph{weak} and \emph{strong survival}.

In \cite{lalley2000} the behaviour of the critical BRW on the free group was
studied in more detail and two phenomena were observed. First,
$\Phi(R)=\mathrm{HD}(\Omega)/2$ if and only if the underlying random walk is
the simple random walk. This statement is not true for our more general setting
since there are non-simple random walks that attain the maximal Hausdorff
dimension $\mathrm{HD}(\Omega)/2$, see Remark \ref{rem:=1/2} together with
Example 3.14. Second, it was shown in \cite{lalley2000} that $\Phi(R)-\Phi(\lambda)\sim C \sqrt{R-\lambda}$ as $\lambda\uparrow R$. For free products of groups this behaviour turns out to be more subtle:  $\Phi(R)-\Phi(\lambda)$ may behave like $C (R-\lambda)$  or $C \sqrt{R-\lambda}$  depending on whether the Green function is differentiable at its radius of convergence or not.

The very same phenomena were also studied in the continuous setting. Lalley and Sellke \cite{lalley:97} studied the phase transition for branching Brownian motion on the hyperbolic disc and Karpelevich,  Pechersky, and Suhov \cite{karpelevich:98} generalized these results to higher  dimensional Lobachevsky spaces. Grigor'yan and Kelbert \cite{grigoryan:03} studied recurrence and transience for branching diffusion processes on Riemannian manifolds. In Cammarota and Orsingher \cite{cammarota} first results on a ``linear'' growing system of particles on the hyperbolic disc are given.

In the case of free products of groups  $\Gamma=\Gamma_1\ast \ldots \ast
\Gamma_r$, where at least one of the factors is infinite, another phase
transition occurs. The boundary $\Omega$ can be decomposed into up to $r+1$
direct summands. For $1\leq i\leq  r$, let  $\Omega_{i}$ denote the set of ends
described by semi-infinite non-backtracking paths, which eventually stay in one
copy of $\Gamma_i$. The set $\Omega_{\infty}$  consists of  all ends described
by infinite, non-backtracking paths that
change the different copies of the free factors infinitely many times. 
Now, for all infinite $\Gamma_{i}$, Theorem \ref{thm:phase-transition} gives a criterion whether $\Lambda \cap \Omega_i\neq \emptyset$ almost surely. In particular, it states that there exists a critical value $\lambda_{i}$
such that $\lambda\leq \lambda_{i}$ is equivalent to $\Lambda \cap \Omega_i =
\emptyset$ almost surely. In other words, if we increase the growth parameter $\lambda$ then more and more different parts of the boundary appear in $\Lambda$.
However, even if $\Lambda\cap\Omega_{i}\neq \emptyset$, only the infinite words contribute to the Hausdorff dimension of $\Lambda$, see Corollary \ref{cor:omega-i-do-not-contribute}.

Finally, for the case of free products of \textit{finite} groups we slightly
adapt the metric defined on the boundary and get (following analogously the reasoning in \cite{lalley2000}) a simpler formula for the Hausdorff dimension of $\Lambda$, see Corollary \ref{Cor:finite-groups}. Analogously, we obtain  a formula for the Hausdorff dimension of $\Lambda$ if we have a BRW on free products by amalgamation of \textit{finite} groups, see Corollary \ref{Cor:amalgam}. In both cases the Hausdorff dimension can be expressed through  a Perron--Frobenius eigenvalue.
\par
Let us remark that free products have been studied in great variety. 
Asymptotic behaviour of return probabilities of random walks on free products has  been
studied in many ways; e.g. Gerl and Woess \cite{gerl-woess}, \cite{woess3}, Sawyer \cite{sawyer78}, Cartwright and Soardi \cite{cartwright-soardi}, \mbox{Lalley \cite{lalley93},}  and Candellero and Gilch \cite{candellero-gilch}. 
For free products of finite groups, Mairesse and Math\'eus \cite{mairesse1} computed an explicit formula for the drift and asymptotic entropy. Gilch \cite{gilch:07}, \cite{gilch:11} computed different formulas for the drift and also for the entropy for random walks on free products of graphs. 
Our proofs envolve in a very crucial way generating functions techniques for free products. These techniques were introduced independently and simultaneously by Cartwright and Soardi \cite{cartwright-soardi}, Woess \cite{woess3}, Voiculescu \cite{voiculescu}, and McLaughlin \cite{mclaughlin}.
In particular, we  show that the Hausdorff dimension can be computed as the solution of a functional equation in terms of double generating functions.
\par
The structure of the paper is as follows. In Section \ref{sec:brw} we  give an introduction to random walks on free products, generating functions, and branching random walks. In Section \ref{sec:results} we state our results and illustrate them with sample computations. The proofs are given in Section \ref{sec:proofs}.

\section{Branching Random Walks on Free Products}
\label{sec:brw}

\subsection{Free Products of Groups and Random Walks}
\label{subsec:fp-rw}
Let $\mathcal{I}=\{1,2,\dots,r\}$ be a finite index set. Suppose we are given
finitely generated groups $\Gamma_{i}$, $i\in\mathcal{I}$, where each $\Gamma_i$ is
generated by a symmetric generating set $S_i$ (that is, $s\in S_i$
implies \mbox{$s^{-1}\in S_i$)} with identity $e_i$. Let $\Gamma_i^\times := \Gamma_i\setminus\{e_i\}$ for every $i\in\mathcal{I}$ and let
$\Gamma_\ast^\times :=\bigcup_{i\in\mathcal{I}}\Gamma_i^\times$. The \textit{free product} $\Gamma:=\Gamma_1\ast \ldots \ast
\Gamma_r$ is defined as the set
\begin{equation}\label{freeproduct}
 \Big\lbrace x_1x_2\dots x_n\ \Bigl|\
n\in\mathbb{N}, x_j\in \Gamma_\ast^\times, x_j\in \Gamma_k^\times \Rightarrow
x_{j+1}\notin \Gamma_k^\times\ \Big\rbrace\cup \Big\lbrace e \Big\rbrace. 
\end{equation}
That is, each element of $\Gamma$ is a \textit{word} $x_1\dots x_n$ such that each
letter (also called \textit{block}) $x_i$ is a non-trivial element of one of the $\Gamma_i$'s and two
consecutive letters are not from the same free factor $\Gamma_i$; $e$ denotes
the empty word.
We exclude the trivial cases where $\Gamma_{i}$ is the trivial group and the case 
$r=2=|\Gamma_1|=|\Gamma_2|$; see beginning of Subsection \ref{subsec:generating-functions} for further remarks.  The group operation on the free product $\Gamma$ can be described as follows: 
if $u=u_1\dots u_m, v=v_1\dots v_n\in \Gamma$ then $uv$ stands for their
concatenation as words with possible contractions and cancellations in the
middle in order to get the form of (\ref{freeproduct}). For instance, if $u=aba$ and $v=abc$ with $a,c\in\Gamma_1^\times, b\in\Gamma_2^\times$ and
$a^2=e_1,b^2\neq e_2$, then $uv=(aba)(abc)=a(b^2)c$.
In particular, we set $ue_i:=u$ for all
$i\in\mathcal{I}$ and $eu:=u$. Note that $\Gamma_i\subseteq \Gamma$ and $e_i$ as a word in
$\Gamma$ is identified with $e$. The \textit{block length} of a word $u=u_1\dots u_m\in\Gamma$
is given by $\Vert u\Vert:=m$.
Additionally, we set $\Vert e\Vert:=0$. The
\textit{type} $\tau(u)$ of $u$ is defined to be $i$ if $u_m\in
\Gamma_i^\times$; we set $\tau(e):=0$.

\begin{figure}
\begin{center}
\includegraphics[width=6cm]{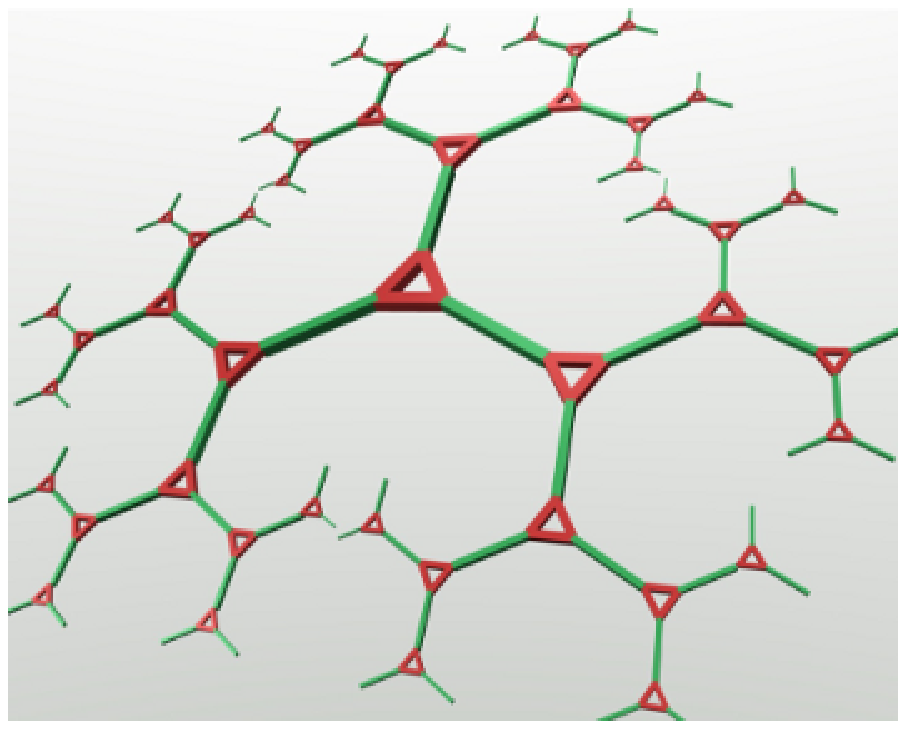}
\end{center}
\caption{Structure of the free product $(\Z/2\Z)\ast(\Z/3\Z)$.}
\label{fig:free-product}
\end{figure}

To help visualizing the structure of a free product we
may interpret the set $\Gamma$ as the vertex set of its Cayley graph $\mathcal{X}$ (with respect to  the generating set $\bigcup_{i\in\calI} S_i$), which is constructed as follows: consider Cayley graphs $\mathcal{X}_1,\dots,\mathcal{X}_r$ of
$\Gamma_1,\dots,\Gamma_r$ w.r.t. the (finite) symmetric  generating
sets $S_1,\dots,S_r$;  take copies of
$\mathcal{X}_1,\dots, \mathcal{X}_r$ and glue them together at their identities to one single common
vertex, which becomes $e$; inductively, at each vertex $v=v_1\dots v_k$ with
$v_k\in \Gamma_i$ attach a copy of every $\mathcal{X}_j$, $j\neq i$, where $v$
is identified with $e_j$ of the new copy of $\mathcal{X}_j$; see Figure
\ref{fig:free-product}. The
natural graph distance on $\mathcal{X}$ is also used for elements of $\Gamma$
and we write $l(u)$ for the  \textit{graph distance} or \textit{length} of
$u\in\Gamma$ to $e$. A \textit{geodesic} of $u$ is a shortest path from $e$ to
$u$. We remark that the length of an element may differ drastically from its
block length. 

\par
We construct in a natural way a random walk on $\Gamma$ from some given random
walks on its free factors. Suppose we are given (symmetric, finitely supported)
probability measures $\mu_i$ on $\Gamma_{i}$ with $\langle
\mathrm{supp}(\mu_i)\rangle=\Gamma_i$ for each $i\in\mathcal{I}$.
For $x,y\in\Gamma_i$, the corresponding single step transition
probabilities of a random walk on $\Gamma_i$ are given by $p_i(x,y):=\mu_i(x^{-1}y)$ and the $n$-step transition
probabilities are denoted by  $p_i^{(n)}(x,y):=\mu^{(n)}_i(x^{-1}y)$, where
$\mu_i^{(n)}$ is the $n$-th convolution power of $\mu_i$. Each of
these random  walks is irreducible.  For sake of simplicity, we  also
assume $\mu_i(e_i)=0$ for every $i\in\mathcal{I}$. We lift $\mu_i$ to a
probability measure $\bar\mu_i$ on $\Gamma$ by defining  $\bar\mu_i(x):=\mu_i(x)$, if $x\in\Gamma_i$, and
$\bar\mu_i(x):=0$ otherwise. Let $\alpha_{i}>0$, $i\in\calI$, with $\sum_{i\in\mathcal{I}} \alpha_i =
1$. We now obtain a new finitely supported probability measure on $\Gamma$ given by
$$
\mu=\sum_{i\in\mathcal{I}} \alpha_i \bar \mu_i.
$$
The random walk on $\Gamma$ starting at $e$, which is governed by $\mu$, is described by the sequence of random variables
$(X_n)_{n\in\mathbb{N}_0}$. For $x,y\in \Gamma$, the associated single and $n$-step transition
probabilities are denoted by $p(x,y):=\mu(x^{-1}y)$ and
$p^{(n)}(x,y):=\mu^{(n)}(x^{-1}y)$, where $\mu^{(n)}$ is the $n$-th convolution
power of $\mu$. The Cayley graph under consideration will always be with respect to the  set of generators $\mathrm{supp}(\mu) = \bigcup_{i\in\mathcal{I}} \mathrm{supp}(\mu_i)$. We refer to Remark \ref{rem:<1/2} for a short discussion for the case of non-nearest neighbour random walks. 

\subsection{Generating Functions}
\label{subsec:generating-functions}
One key ingredient of the proofs is the study of  the following generating functions.
The most common among these generating functions are the \textit{Green functions} related to $\mu_i$ and $\mu$ which are defined by $$
G_i(x_i,y_i|z):= \sum_{n\geq 0} p_i^{(n)}(x_i,y_i)\,z^n \quad \textrm{ and } \quad 
G(x,y|z):= \sum_{n\geq 0} p^{(n)}(x,y)\,z^n,
$$
where $z\in\mathbb{C}$, $i\in\mathcal{I}$, $x_i,y_i\in \Gamma_i$ and $x,y\in \Gamma$. We note that the free product $\Gamma$ is non-amenable and that 
the radius of convergence $R$ of
$G(\cdot,\cdot|z)$ is strictly larger than $1$; see e.g. \cite[Thm. 10.10, Cor. 12.5]{woess}. In particular, this implies \textit{transience} of our
random walk on $\Gamma$. At this point let us remark that the case
$r=2=|\Gamma_1|=|\Gamma_2|$ leads to a recurrent random walk (and therefore to
a recurrent branching random walk), which is the reason why we excluded this case. Moreover, non-amenability of $\Gamma$ yields $G(e,e|R)<\infty$; see e.g. \cite[Proposition 2.1]{lalley2}.
\par
The \textit{first visit
  generating functions} related to $\mu_i$ and $\mu$ are given by
\begin{eqnarray*}
F_i(x_i,y_i|z) & := & \sum_{n\geq 0} \Prob\bigl[Y_n^{(i)}=y_i,\forall m\leq n-1:
Y_m^{(i)}\neq y_i \mid Y_0^{(i)}=x_i\bigr]z^n\ \textrm{ and}\\
F(x,y|z) & := & \sum_{n\geq 0} \Prob\bigl[X_n=y,\forall m\leq n-1:
X_m\neq y \mid X_0=x\bigr]z^n,
\end{eqnarray*}
where $\bigl(Y_n^{(i)}\bigr)_{n\in\mathbb{N}_0}$ describes a random walk on
$\Gamma_i$ governed by $\mu_i$. For $M\subseteq \Gamma$, we also define
$$
F(x,M|z)  :=  \sum_{n\geq 0} \Prob\bigl[X_n\in M,\forall m\leq n-1:
X_m\notin M \mid X_0=x\bigr]z^n
$$
and the \textit{first return generating function}
$$
U(x,M|z)  :=  \sum_{n\geq 1} \Prob\bigl[X_n\in M,\forall~1\leq m\leq n-1:
X_m\notin M \mid X_0=x\bigr]z^n.
$$
By a Harnack-type inequality the generating functions $F(\cdot,\cdot|z)$ and $U(\cdot,\cdot|z)$
have also radii of convergence of at least $R>1$ and $U(x,M|z)=F(x,M|z)$ if $x\notin M$. By transitivity, we have $G_i(x_i,x_i|z)=G_i(e_i,e_i|z)$ and $G(x,x|z)=G(e,e|z)$ for
all $x_i\in\Gamma_i$ and $x\in\Gamma$. 
For $x\in\Gamma\setminus\{e\}$, we have 
\begin{equation}\label{equ:F<1}
G(e,e|z)>F(e,x|z) G(x,x|z) F(x,e|z);
\end{equation}
indeed, while on the left hand side we take into account all paths from $e$ to
$e$, on the right hand side we only take into account all random walk paths
from $e$ to $e$ which pass through $x$; therefore, strict inequality
follows from irreducibility of the random walk which ensures always existence
of random walk paths from $e$ to $e$ not passing through $x$ . 
Symmetry of the $\mu_i$'s now implies that $F(e,x|z)<1$
for all $|z|\leq R$ and all $x\in\Gamma\setminus \{e\}$. The \textit{last visit
  generating functions} related to $\mu_i$ and $\mu$ are given by
\begin{eqnarray*}
L_i(x_i,y_i|z) & := & \sum_{n\geq 0} \Prob\bigl[Y_n^{(i)}=y_i,\forall 1\leq
m\leq n:
Y_m^{(i)}\neq x_i \mid Y_0^{(i)}=x_i\bigr]z^n\ \textrm{ and}\\
L(x,y|z) & := & \sum_{n\geq 0} \Prob\bigl[X_n=y,\forall 1\leq m\leq n:
X_m\neq x \mid X_0=x\bigr]z^n.
\end{eqnarray*}
We have the following important equations, which follow by conditioning on 
the first visits of $y_i$ and $y$, the last visits of $x_i$ and $x$ respectively: 
\begin{equation}\label{gl-equations}
\begin{array}{rcl}
G_i(x_i,y_i|z)& = & F_i(x_i,y_i|z)\cdot G_i(y_i,y_i|z)=G_i(x_i,x_i|z)\cdot L_i(x_i,y_i|z),\\[1ex]
G(x,y|z) & = & F(x,y|z)\cdot G(y,y|z)=G(x,x|z)\cdot L(x,y|z).
\end{array}
\end{equation}
Thus, by transitivity we obtain 
\begin{equation}
F(x,y|z)=L(x,y|z) \mbox{ for any } x,y\in\Gamma \mbox{ and } |z|\leq R.
\end{equation}
Let $x,y,w\in \Gamma$ such that all (random walk) paths from $x$ to $w$ pass through $y$. Then 
\begin{equation}\label{equ:f-l-product}
F(x,w|z) = F(x,y|z) \cdot F(y,w|z) \quad \textrm{ and }\quad L(x,w|z) = L(x,y|z) \cdot L(y,w|z);
\end{equation}
this can be checked by conditioning on the first/last visit of $y$ when
walking from $x$ to $w$.
For $i\in\calI$ and $z\in\mathbb{C}$, we define the functions
\begin{equation}
\xi_i(z) :=  U\bigl(e,\mathrm{supp}(\mu_i)|z\bigr) = U\bigl(e,\Gamma_{i}^\times|z\bigr) = F\bigl(e,\mathrm{supp}(\mu_i)|z\bigr),
\end{equation} 
which have also radii of convergence of at least $R>1$. We remark that $\xi_i(1)<1$; see e.g. \cite[Lemma 2.3]{gilch:07}. Moreover, we have
$F(x_i,y_i|z)=F_i\bigl(x_i,y_i|\xi_i(z)\bigr)$ and $L(x_i,y_i|z)=L_i\bigl(x_i,y_i|\xi_i(z)\bigr)$ for all $x_i,y_i\in
\Gamma_i$; see \cite[Prop. 9.18c]{woess} and \cite[Lemma 2.2]{gilch:07}. Thus, by conditioning on the number of visits of $e$ before finally making a step from $e$ to $\Gamma_i^\times$ we get the following formula:
\begin{equation}\label{equ:xi-formula}
\xi_i(z) = \frac{\alpha_i z}{1-\sum_{j\in\mathcal{I}\setminus\{i\}}
  \sum_{s\in \Gamma_j} \alpha_j \mu_j(s)z F_j\bigl(s,e_j \bigl| \xi_j(z)\bigr)}.
\end{equation}

Finally, we define the following power series that will lead to a useful expression for the Hausdorff dimension. Let
\begin{equation}\label{equ:scriptF}
\mathcal{F}(\lambda|z):=\sum_{x\in\Gamma} F(e,x|\lambda)\,z^{l(x)},
\end{equation} and define for $i\in\mathcal{I}$:
\begin{eqnarray}
\mathcal{F}_i^+(\lambda|z) & := & \sum_{x\in\Gamma_i^\times} F(e,x|\lambda)\,z^{l(x)}=
\sum_{x\in\Gamma_i^\times} F_i\bigl(e_i,x\bigl|\xi_i(\lambda)\bigr)\,z^{l(x)},\label{equ:scriptFi+}\\
\mathcal{F}_i(\lambda|z) & := & \sum_{n\geq 1}\sum_{\substack{x=x_1\dots x_n
  \in\Gamma:\\ x_1\in\Gamma_i^\times}} F(e,x|\lambda)\,z^{l(x)}
= \mathcal{F}_i^+(\lambda|z)\Bigl(1+\sum_{j\in\mathcal{I}\setminus\{i\}}
\mathcal{F}_j(\lambda|z)\Bigr).
\label{equ:scriptFi}
\end{eqnarray}
The latter functions satisfy the following relation:
\begin{equation}\label{equ:scriptF2}
\mathcal{F}(\lambda|z) = 1 + \sum_{i\in\calI} \mathcal{F}_i(\lambda|z).
\end{equation}

\subsection{Branching Random Walks}
In this subsection we introduce discrete-time branching random walks on free
products and recall some basic results. 

There are two different main descriptions or constructions of a branching random walk (BRW). The first defines the process inductively  as a \textit{growing cloud} of particles moving in (discrete) time and space. The second, via tree-indexed random walks, uses the fact that the branching distribution does not depend on the space. For that reason one can separate branching and movement into two steps. First, one generates the whole genealogy of the process and then one maps the corresponding genealogical tree into the Cayley graph. 
In both cases we need the following definition.  A Galton--Watson  process  is
characterized through an \textit{offspring distribution} $\nu$. This is a
probability measure  on $\mathbb{N}=\{0,1,2,3,\dots\}$ with mean (or also
called growth parameter) $\lambda=\sum_{k\geq 1} k\,\nu(k)\in (0,\infty)$. We
assume that $\nu$ has \textit{finite second moment}, that is, $\sum_{k\geq 1}
k^2\,\nu(k)<\infty$.  Moreover we exclude the cases where $\nu(0)>0$ and
$\nu(1)=1$; this  guarantees that the process survives almost surely and that the BRW is not reduced to a (non-branching) random walk. 

The BRW on $\Gamma$ is defined inductively: at time $0$ we have one particle at
$e$ (if not mentioned otherwise). Between time $n$ and $n+1$ the process performs two steps: branching and movement. First,  each particle,  independently of all others and the previous history of the process,  produces descendants according to $\nu$ and dies.  Second, each of these descendants, independently of all others and the past, moves  to a neighbour vertex in $\Gamma$ according to $\mu$. 
A particle located at some vertex $x\in\Gamma$ at time $n$ has a unique direct
ancestor at time $n-1$. Consequently, each particle  has a unique finite sequence of ancestors, the family history, which traces back to the original starting particle at $e$. The
sequence of the locations of its ancestors (chronologically ordered) gives a
path from $e$ to $x$, which we call the \textit{trail} of the particle. 
\par
Sometimes it will be convenient  to work with the interpretation of a BRW as a tree-indexed random walk, see \cite{benjamini:94}.
Let $\T$ be a rooted infinite tree. The root is denoted by $\ro$ and other vertices by $v$ and let $|v|$ be the (graph) distance from $v$ to the root $\ro$. 
The random walk on $\Gamma$ indexed by $\T$ is the collection of
$\Gamma$-valued random variables $(S_{v})_{v\in \T}$ defined as follows. Label
the edges of $\T$ with i.i.d. random variables $\eta_v$ with distribution
$\mu$; the random variable $\eta_v$ is the label of the edge $(v^-,v)$. Define
$S_v=e \cdot \prod_{i=1}^{|v|} \eta_{v_i}$ where $\langle v_0=\ro, v_1, \ldots,
v_n=v\rangle$ is the unique geodesic (also called ancestry of $v$) from $\ro $ to $v$ at level $n$.
A tree-indexed random walk becomes a BRW if the underlying tree is a
Galton--Watson tree induced by $\nu$.  
We refer to  $\T$ as the family tree  and to $\mathcal{X}$  as the base graph of the BRW. Furthermore, a vertex $v\in \T$ is called a particle of the BRW and $\T_{n}$ denotes  the vertices of $\T$ on level $n$ or equivalently the particles in generation $n$.

A useful variation of the first description of a BRW is the  \textit{coloured
  branching random walk}, see \cite{lalley2000}. This process behaves like a
standard BRW where in addition each particle is either blue or red. In order to
define this coloured version we choose a subset $M$ of $\Gamma$ that plays the
role of a ``paint bucket''. We start the BRW with one blue particle at
$e$. Blue particles located outside of $M$ produce blue offspring. A blue particle that hits the paint bucket is frozen there and will be replaced by a red particle. The new red particle starts an ordinary (red-coloured) branching random walk.  As a consequence, every red particle has exactly one ``frozen'' ancestor in $M$.

We denote by $Z_\infty(M)\in \N\cup \{\infty\}$ the random number of frozen (blue) particles in $M$ during the whole branching process. If $M=\{x\}$ then we just write $Z_\infty(x)$.
\par
For ease of presentation we will switch freely between the different definitions of a BRW; nevertheless it will  always be clear from the context which description we are using.
\par
A BRW on a Cayley graph is called \emph{recurrent} if each vertex is visited
infinitely many times and \emph{transient} if any finite subset is eventually
free of particles. The recurrence/transience behaviour is well understood. In
fact,  we have the following classification in recurrence and transience, see
\cite {benjamini:94} for the sub- and supercritical case and \cite{gantert:04} for the critical case. We also refer to \cite{grigoryan:03} for the corresponding result in the continuous setting.
\begin{Th}
The BRW is transient if and only if $\lambda\leq R$.
\end{Th}
Recall that in the language of \cite{lalley2000} transience is equivalent to
\emph{weak survival} if $\lambda >1$. For the rest of this paper we will restrict our investigation to the case of transience or weak survival. Since in this case the process eventually vacates every finite subset of
$\Gamma$  almost surely  the investigation of the \emph{convergence} of the BRW to the geometric boundary is meaningful.

\subsection{Ends of Graphs, Box-Counting Dimension and Hausdorff Dimension}\label{subsec:endofgraphs}
Let us first recall some basic notations on infinite graphs. Let $\mathcal{G}$ be an infinite, connected, locally finite graph  with countable vertex set and root $e$. For ease of presentation, we  will identify $\mathcal{G}$ or a subgraph with its vertex set. A
\textit{path} of length $n$ in $\mathcal{G}$ is a finite sequence of vertices $[x_0,x_1,\dots,x_n]$
such that there is an edge from $x_{i-1}$ to $x_i$ for each
$i\in\{1,\dots,n\}$. Recall that a \textit{geodesic} of a vertex $x\in \mathcal{G}$ is a shortest path
from $e$ to $x$ in $\mathcal{G}$.
A \textit{ray} is a semi-infinite path
$[e=x_0,x_1,x_2,\dots]$, which does not backtrack, that is, $x_i\neq x_j$ if
$i\neq j$. 
Two rays $\eta_1$ and $\eta_2$ are \textit{equivalent} if there is a
third ray which shares infinitely many vertices with  $\eta_1$ and $\eta_2$. An equivalence
class of rays is called an \textit{end}. The set of
equivalence classes of rays is called the \textit{end boundary} of $\mathcal{G}$, denoted
by $\partial \mathcal{G}$. For further details we refer to \mbox{\cite[Section 21]{woess}.}
\par
In the case of free products we have different types of ends occurring in the Cayley graph $\mathcal{X}$ of $\Gamma$: ends arising from ends in one of the
$\mathcal{X}_i$, and ``infinite words''. More precisely, denote by $\Omega_i^{(0)}$ the set of ends
of $\mathcal{X}_i$. For $\omega_i\in\Omega_i^{(0)}$, let 
$\eta=[e_i,y_1,y_2,\dots]\in \omega_{i}$ and let  $x\in\Gamma$,
where $[x_0,x_1,\dots,x_n]$ is a geodesic from $x_{0}$ to $x=x_{n}$. Then, the ray
$x\eta:=[x_0,x_1,\dots,x_n,x_ny_1,x_ny_2,\dots]$ describes an end in
$\Gamma$. The end described by $x\eta$ is denoted by $x\omega_i$. We set $\Omega_i:=\bigl\lbrace x\omega_i \mid x\in\Gamma,\omega_i\in\Omega_i^{(0)} \bigr\rbrace$. Moreover, the set of infinite words is given by
$$
\Omega_\infty = \bigl\lbrace
x_1x_2x_3\ldots \in\bigl(\Gamma_\ast^\times\bigr)^{\mathbb{N}}\,\bigl|\,
x_j\in\Gamma_k^\times \Rightarrow x_{j+1}\notin\Gamma_k^\times
\bigr\rbrace.
$$
It is easy to see that the set $\Omega$ of ends of $\mathcal{X}$ can be decomposed
in the following way:
$$
\Omega = \Omega_\infty \uplus \Omega_1 \uplus \Omega_2 \uplus \dots \uplus \Omega_r.
$$
Observe that $\Omega_i$ is empty if and only if $\Gamma_i$ is finite. Thus, if all groups
$\Gamma_i$ are finite then $\Omega=\Omega_\infty$.
\par
In order to measure the size of $\Omega$ we  define a metric on
$\Omega$.   We say that an  end $\omega_{1}\in \Omega$ is contained in a subset of the graph  if all representatives have all but finitely many vertices in this subset. Now, if we remove from $\mathcal{X}$ any finite vertex subset $F\subseteq
\mathcal{X}$ (including the removal of edges to vertices in $F$) then there is exactly
one connected component in the reduced graph $\mathcal{X}\setminus F$ containing the end $\omega_1$. We call this component the $\omega_1$-component and say that $\omega_{1}$ ends up in this component.
 Denote by $B_m:=\{x\in\Gamma \mid l(x)\leq m\}$ the ball
centered at $e$ with radius $m$; we also set $B_{-1}:=\emptyset$. Let $\omega_2\in\Omega$ be another end with
$\omega_1\neq \omega_2$. Obviously, there is some maximal $m\in\N_0$ such that
$\omega_1$ and $\omega_2$ end up in the same connected component of $\mathcal{X}\setminus B_{m-1}$. We write $c(\omega_1,\omega_2)$ for this maximal integer $m$.
We now define a
metric on $\Omega$ by
$$
d_\Omega(\omega_1,\omega_2):= \alpha^{c(\omega_1,\omega_2)},
$$
where $\alpha\in(0,1)$ is arbitrary, but fixed. Additionally, we set $d_\Omega(\omega_1,\omega_1):=0$. The ball $B(\omega,\eps)$ centered
at $\omega\in\Omega$ with radius $\eps\geq 0$ is given by all ends $\hat\omega\in\Omega$
with $d_\Omega(\omega,\hat\omega)\leq \eps$. In other words, if $\eps=\alpha^m$ then
$\hat\omega\in B(\omega,\eps)$ if and only if $\omega$ and $\hat\omega$ end up in
the same component of $\mathcal{X}\setminus  B_{m-1}$. 
\par
A \textit{cover} of a subset $\Omega'\subseteq \Omega$ is a finite or
countable set of balls of the form $B(\omega,\eps_\omega)$ with $\omega\in\Omega'$ and
$\eps_\omega>0$ such that the union of these balls include $\Omega'$.
 For any $\eps>0$ let $N_\eps(\Omega')$ be
the minimal number of balls of the form $B(\omega,\eps_\omega)$ with
$\omega\in\Omega'$ and $0<\eps_\omega\leq \eps$, which cover $\Omega'$. Apparently,
$N_\eps(\Omega')$ is bounded from above by the number of elements in $\Gamma$ at
graph distance $m=\lceil \log(\eps)/\log(\alpha)\rceil$. The
\textit{lower} and \textit{upper box-counting dimension} (also called \textit{Minkowski dimension})
of $\Omega'$ are defined as
\begin{equation}\label{def:BD}
\mathrm{\underline{BD}}(\Omega') := \liminf_{\eps\downarrow 0} \frac{\log
  N_\eps(\Omega')}{-\log \eps} \quad \mbox{ and } \quad \mathrm{\overline{BD}}(\Omega') := \limsup_{\eps\downarrow 0} \frac{\log N_\eps(\Omega')}{-\log \eps}. 
\end{equation}
If both limits are equal then the common value is called the
\textit{box-counting dimension} $\mathrm{BD}(\Omega')$ of $\Omega'$.

Another well-known measure for the size of $\Omega'$ is given by the
Hausdorff dimension. For $\delta>0$, the \textit{$\delta$-dimensional
Hausdorff measure} of $\Omega'$ is defined by
$$
\mathcal{H}_\delta(\Omega') := \lim_{\varepsilon\downarrow 0}\inf\Bigl\lbrace \sum_i \eps_i^\delta\,\Bigl|\, \bigl\{B(\cdot,\eps_i)\bigr\}_i \textrm{ is a cover of } \Omega' \textrm{ with } \eps_i<\varepsilon \Bigr\rbrace.
$$
Then the \textit{Hausdorff dimension} of $\Omega'$ is defined as
\begin{equation}\label{def:HD}
\mathrm{HD}(\Omega') := \inf \bigl\lbrace
\delta \geq 0 \bigl|\, \mathcal{H}_\delta(\Omega') =0 \bigr\rbrace.
\end{equation}
Since $\mathcal{X}$ has bounded vertex degrees we have $\mathrm{HD}(\Omega')<\infty$. It is well-known that, for all $\Omega'\subseteq \Omega$,
$$
\mathrm{HD}(\Omega')\leq \underline{\mathrm{BD}}(\Omega').
$$
One of our main goals is to investigate to which kind of ends the branching random walk
converges and to compare  the dimensions of the whole space of ends with the
set of ends which are ``hit'' by the BRW. More precisely, for any $\omega\in \Omega$, if we remove any
finite vertex subset $F\subseteq \mathcal{X}$ then there is exactly
one connected component in the reduced graph $\mathcal{X}\setminus F$ containing $\omega$.
We say that the branching
random walk \textit{accumulates} at the end $\omega$ if for every finite vertex
subset $F\subseteq \mathcal{X}$ there is at least one particle visiting the connected
$\omega$-component in $\mathcal{X}\setminus F$. The set of accumulation points is denoted
by $\Lambda$. If the BRW is recurrent then $\Omega=\Lambda$;
thus, we restrict our investigation to the more interesting case of transience and therefore assume  $1<\lambda\leq R$. Note that $\Lambda\cap\Omega_\infty$ is almost surely non-empty;  each infinite ancestry line  converges to some element in $g_1g_2\dots\in\Omega_\infty$ with convergence in the sense that the
length of the common prefix of the particle's location and $g_1g_2\dots$ tends to
infinity, see e.g. \cite[Proposition 2.5]{gilch:07}.  We remark also that the Hausdorff dimensions of $\Lambda$ and $\Lambda\cap \Omega_\infty$ are almost surely constant, which can be shown analogously as explained in \cite[Sec.~1, Remark (C)]{lalley2000}.

\section{Results}\label{sec:results}

In this section we summarize our results about branching random walks on free products and present several explicit examples.

\subsection{Main Results}\label{subsec:main-results}
The  first result describes how the structure of $\Lambda$ gets richer when
increasing the growth parameter $\lambda$ and that there are up to $r=|\calI|$ possible phase transitions.
\begin{Th}\label{thm:phase-transition}
Let $\lambda\in(1,R]$. Then $\Prob\bigl[\Lambda\cap \Omega_i\neq
\emptyset\bigr]\in\{0,1\}$, and $\Prob\bigl[\Lambda\cap \Omega_i\neq
\emptyset\bigr]=1$ if and only if $\xi_i(\lambda)> 1$. More precisely:
\begin{enumerate}
\item If $\xi_i(\lambda)\leq 1$ then $\emptyset\subsetneq \Lambda\subseteq \Omega_\infty$.
\item If $\xi_i(\lambda)> 1$ then $\emptyset\subsetneq\Omega_\infty \cap \Lambda \subset \Lambda$ with $\Lambda
  \cap \Omega_i\neq \emptyset$ and $|\Lambda \cap \Omega_i|=\infty$.
\end{enumerate}
\end{Th}

\begin{Rem}
In the case where one of the free factors is an infinite amenable group its ends do not appear in $\Lambda$. In other words, if $R_i=1$ is the radius of convergence of $G_i(e_i,e_i|z)$ then $\xi_i(\lambda)\leq 1$ for all $\lambda\in(1,R]$; see \cite[Lemma
  17.1a]{woess}. Consequently, no ends in $\Omega_i$ contribute to $\Lambda$,
  that is, $\Lambda\cap\Omega_i=\emptyset$ almost surely.
\end{Rem}

We illustrate the above described behaviour in the following two examples:
\par
\addtocounter{Th}{1}
\begin{Example}
Consider  $\Gamma=\Z^{d_1}\ast \Z^{d_2}$ and let   $\mu_1$ and $\mu_2$ be two symmetric probability measures on $\Z^{d_1}$ and $\Z^{d_2}$. Due to Kesten's amenability criterion we have $R_{1}=R_{2}=1$.  Consequently, $\Lambda \subseteq \Omega_\infty$ almost surely for all $\lambda \leq R$.  
 \end{Example}
 \addtocounter{Th}{1}
\begin{Example}
Consider $\Gamma=\Gamma_1\ast\Gamma_2$, where $\Gamma_{1}$ and $\Gamma_{2}$ are non-amenable groups, and let $\mu_{i}$ define a symmetric random walk on $\Gamma_{i}$ for $i\in\{1,2\}$. Due to the non-amenability we have that  $R_1,R_2>1$ and $G_i(e_i,e_i|R_i)<\infty$.  In the case where 
$$
\alpha_1=\frac{R_1\,G_1(e_1,e_1|R_1)}{R_1\,G_1(e_1,e_1|R_1)+R_2\,G_2(e_2,e_2|R_2)}
$$ 
we obtain by \cite[Lemma 17.1]{woess} that $\xi_1(R),\xi_2(R)>1$. Therefore, there are numbers \mbox{$\lambda_1,\lambda_2\in(1,R)$} with $\xi_1(\lambda_1)=\xi_2(\lambda_2)=1$ which leads to phase transitions at $\lambda_{1}$ and $\lambda_{2}$. 
\end{Example}
Now we state our first main result.
\begin{Th}\label{thm:boxdim}
Suppose that $\nu$ has finite second moment. Then the box-counting dimension of $\Lambda$, $\Lambda\cap\Omega_\infty$
respectively, exists and  equals the Hausdorff dimension of $\Lambda$, $\Lambda\cap\Omega_\infty$
respectively. Furthermore:
$$
\mathrm{BD}(\Lambda)=\mathrm{BD}(\Lambda \cap \Omega_\infty)=\mathrm{HD}(\Lambda)=\mathrm{HD}(\Lambda \cap \Omega_\infty)=\frac{\log z^\ast}{\log\alpha},
$$
where $z^\ast$  is the smallest real positive number
 with 
\begin{equation}\label{equ:z-ast}
\sum_{i\in\mathcal{I}}
\frac{\mathcal{F}_i^+(\lambda|z^\ast)}{1+\mathcal{F}_i^+(\lambda|z^\ast)}=1.
\end{equation}
\end{Th}
\begin{Rem}
The proof of Theorem \ref{thm:boxdim} directly applies to BRW on free products
of \textit{finite graphs} and a corresponding result holds verbatim; see
e.g. \cite[Sec. 9.C]{woess} for a formal definition of general free products and random
walks on them.
\end{Rem}
As a first consequence we obtain that only infinite words contribute to the dimension of $\Lambda.$
\begin{Cor}\label{cor:omega-i-do-not-contribute}
For $i\in\calI$, $\mathrm{HD}(\Lambda \cap \Omega_i) < \mathrm{HD}(\Lambda \cap \Omega_\infty)$. 
\end{Cor}
For $i\in\calI$, $m\in\N$ and $z\in\mathbb{C}$, we define $S_i(m):=|\{x\in\Gamma_i \mid l(x)=m\}|$ and
$$
\mathcal{S}_i^+(z):=\sum_{m\geq 1} S_i(m) z^m.
$$
Analogously to Theorem \ref{thm:boxdim}, we can prove existence of the box-counting dimension of the whole boundary $\Omega$ and express the dimension as the solution of a functional equation.
\begin{Th}\label{thm:boxdim-omega}
The box-counting dimensions of $\Omega$ and $\Omega_\infty$ exist and satisfy
$$
\mathrm{BD}(\Omega)=\mathrm{BD}(\Omega_\infty)=\mathrm{HD}(\Omega)=\mathrm{HD}(\Omega_\infty)=
\frac{\log z^\ast_{\mathcal{S}}}{\log\alpha},
$$
where $z^\ast_{\mathcal{S}}$ is the smallest real positive number with 
\begin{equation}\label{equ:z-ast-S}
\sum_{i\in\mathcal{I}}
\frac{\mathcal{S}_i^+(z^\ast_{\mathcal{S}})}{1+\mathcal{S}_i^+(z^\ast_{\mathcal{S}})}=1.
\end{equation}
\end{Th}
Analogously to Corollary \ref{cor:omega-i-do-not-contribute} we obtain that the Hausdorff dimension of $\Omega$ arises only from the ends in $\Omega_\infty$.
\begin{Cor}\label{cor:boxdim-omega}
For all $i\in\calI$, $\mathrm{HD}(\Omega_i)<\mathrm{HD}(\Omega_\infty)$.
\end{Cor}
Beyond these first consequences of Theorems \ref{thm:boxdim} and \ref{thm:boxdim-omega}, the expressions for the Hausdorff dimensions allow us to study first regularity properties. For any fixed free product $\Gamma$, let us consider the function 
$$
\Phi: [1,\infty)\to \mathbb{R}: \lambda \mapsto \mathrm{HD}(\Lambda),
$$
which assigns to every value $\lambda$ the Hausdorff dimension of $\Lambda$ of
a BRW with growth parameter $\lambda$. The limit case $\lambda=1$ corresponds
to the degenerate case of a non-branching random walk; in this case the Hausdorff dimension is just zero.
\begin{Th}\label{Th:lambda-function}
The function $\Phi(\lambda)$ has the following properties:
\begin{enumerate}
\item $\Phi(\lambda)$ is strictly increasing on $[1,R]$,
  $\Phi(1)=0$ and $\Phi(\lambda)=\mathrm{HD}(\Omega)$ for all $\lambda >R$.
\item $\Phi(\lambda)$ is continuous in $[1,\infty)\setminus \{R\}$ and continuous from the left at $\lambda=R$ with 
$$
\Phi(R)\leq \frac12\mathrm{HD}(\Omega).
$$
\item $\Phi(\lambda)$ has the following behaviour as $\lambda\uparrow R$:
$$
\Phi(R)-\Phi(\lambda) \sim
\begin{cases}
C_1\cdot (R-\lambda), & \textrm{if } G'(R)<\infty,\\
C_2\cdot \sqrt{R-\lambda}, &  \textrm{if } G'(R)=\infty
\end{cases}
$$
for a suitable constant $C_1$, $C_2$ respectively.
\end{enumerate}
\end{Th}

\begin{Rem}\label{rem:<1/2}
The last theorem states that $\mathrm{HD}(\Lambda)$ does not exceed $\mathrm{HD}(\Omega)/2$ unless the BRW is recurrent. We always assumed the random walk to be of nearest neighbour type. However, we feel confident that our techniques work well in the case of finite range random walks and that the equality  $\mathrm{HD}(\Lambda)\leq\mathrm{HD}(\Omega)/2$ does not depend on the choice of the metric. This type of phenomenon was already conjectured for the contact process on the homogeneous tree  in \cite{lalley:96}. We also refer to Section $8$ in \cite{lalley:97} for a discussion how the value $1/2$ can be explained through the ``backscattering principle''.
\end{Rem}

\begin{Rem}\label{rem:=1/2}
In \cite{lalley2000} it was shown that
$\mathrm{HD}(\Lambda)=\mathrm{HD}(\Omega)/2$ only if $\lambda=R$ and if the
underlying walk is a simple random walk. In our more general setting this is no
longer true, since the maximal Hausdorff dimension can also be attained by a
non-simple random walk, see Example 3.14.
More generally, we conjecture that one has maximal
dimension for the BRW (with $\lambda$ being the critical growth value) for every choice of $\alpha_1\in(0,1)$ if we consider a
general free product $\Gamma=\Gamma_1\ast\Gamma_2$ with $\mu_1$ and $\mu_2$ governing positive recurrent random walks on the single factors $\Gamma_1$ and $\Gamma_2$.
\end{Rem}


\begin{Rem} Recall that we always assume that the random walk on $\Gamma$ is symmetric. This assumption  can be dropped for free products of \textit{finite} groups/graphs. In this case we always have the crucial property $F(e,x|R)<1$ for all $x\in\Gamma\setminus\{e\}$ (compare with (\ref{equ:F<1})). In fact, if $x=x_1\dots x_m\in\Gamma\setminus \{e\}$ then
$$
F(e,x_1\dots x_m|R) = \prod_{j=1}^m F_{\tau(x_j)}\bigl(e_{\tau(x_j)},x_j\mid \xi_{\tau(x_j)}(R) \bigr)<1,
$$
as $\xi_i(R)<1$ due to \cite[Lemma 17.1, Theorem 9.22]{woess}.
\end{Rem}

%
%
Theorem \ref{thm:boxdim}  allows explicit calculations in all cases where formulas for the involved generating functions are known.  In the following examples we  set the exponent of the metric on $\Omega$ equal to $1/2$, i.e., $d_{\Omega}(\cdot,\cdot)=2^{-c(\cdot,\cdot)}$.
%
%
\par
\addtocounter{Th}{1}
\newcounter{helpcounter1}
\newcounter{helpcounter2}
\begin{Example} \setcounter{helpcounter1}{\arabic{section}}\setcounter{helpcounter2}{\arabic{Th}}\label{ex:Z2-Z3}
Consider the free product $\Gamma=\Gamma_1\ast\Gamma_2=(\Z/3\Z)\ast(\Z/2\Z)$,
where $\Z/3\Z=\{e_1,a,a^2\}$, with $\mathrm{supp}(\mu_1)=\{a,a^2\}$. The required generating
functions $F(e,x|\lambda)$, $x\in \Gamma_\ast^\times$, may e.g. be obtained by solving
the finite systems of equations given in \cite[Prop. 3c]{woess3}, and therefore
$\mathrm{HD}(\Lambda)$ can be computed via Equation (\ref{equ:z-ast}). Solving
Equation (\ref{equ:z-ast-S}) leads to $\mathrm{HD}(\Omega)=1/2$. Figure \ref{fig:hd-z2-z3} shows
-- with the help of numerical computation by \textsc{Mathematica} -- the graph
of the function $\lambda\mapsto \mathrm{HD}(\Lambda)$ for simple
random walk on $(\Z/3\Z)\ast(\Z/2\Z)$. Let us remark that in this case the
critical parameter $R$ can be explicitely calculated by the formula given in
\cite[(9.29),(3)]{woess}. 
%
%
%
%
%

\par
Another interesting phenomenon occurs in this example. If
$\mu_1(a)=\mu_1(a^2)=1/2$ and if we let $\alpha_1$ vary in the interval $(0,1)$
and denote by
$R(\alpha_1)$ the radius of convergence of $G(e,e|z)$ in dependence of
$\alpha_1$ then we always
get $\Phi\bigl(R(\alpha_1)\bigr)=\frac{1}{2}\mathrm{HD}(\Omega)$, which can be
verified by explicit calculations with the help of \textsc{Mathematica}. 
\begin{figure}
\begin{center}
\includegraphics[width=6cm]{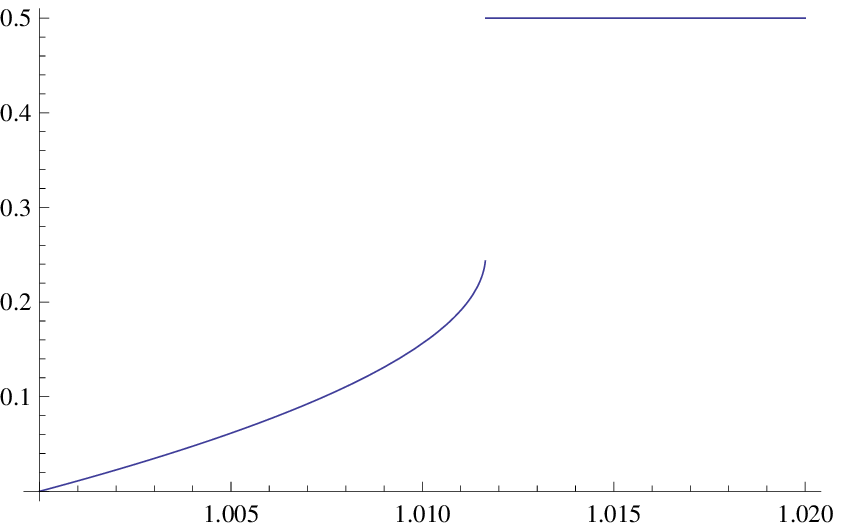}
\end{center}
\caption{Hausdorff dimension $\mathrm{HD}(\Lambda)$ of a  BRW on $(\Z/3\Z)\ast(\Z/2\Z)$ in dependence of $\lambda$ on the $x$-axis}
\label{fig:hd-z2-z3}
\end{figure}
\end{Example}
%
%
%
%
%
\par
\addtocounter{Th}{1}
\begin{Example}
We consider the free product of two infinite ``ladders'' $\Z\times
(\Z/2\Z)$. We set  $\alpha_1=\alpha_2=1/2$ and $\mu_1\bigl((\pm 1,0)
\bigr)=\mu_1\bigl((0,1) \bigr)=\mu_2\bigl((\pm 1,0) \bigr)=\mu_2\bigl((0,1)
\bigr)=1/3$. The functions $F_1\bigl((0,0),(z,a)|z \bigr)$  with
$(z,a)\in\Z\times \Z/2\Z$ can be computed by solving a system of equations as
it is shown in \cite[Section 7.2]{gilch:11}. In order  to compute the Hausdorff
dimension of $\Lambda$ one has to solve,  analogously to \cite[Section 6.2]{gilch:07}:
\begin{eqnarray*}
&&\frac{\lambda}{2} \frac{\xi_1(\lambda)}{\lambda - \xi_1(\lambda)} \\
&=& \frac{\xi_1(\lambda)}{1-\frac{2\xi_1(\lambda)}{3}
  \bigl(F_1\bigl((0,0),(1,0)|\xi_1(\lambda) \bigr)+
  F_1\bigl((0,0),(-1,0)|\xi_1(\lambda)
  \bigr)+F_1\bigl((0,0),(0,1)|\xi_1(\lambda) \bigr)\bigr)}.
\end{eqnarray*}
In order to compute $\mathrm{HD}(\Omega)$ we observe that $S_1(1)=3$ and $S_1(m)=4$ for $m\geq 2$. Hence, $\mathcal{S}_1^+(z)=\mathcal{S}_2^+(z)=3z+4z^2/(1-z)$. This yields $z^\ast_\mathcal{S}=\sqrt{5}-2$. Numerical evaluations then lead to a picture qualitatively similar to Figure  \ref{fig:hd-z2-z3}.
\end{Example}

\subsection{Free Products of Finite Groups}

In this subsection we give a more explicit formula for the box-counting
dimension with respect to  a slightly changed metric on the boundary in the
case of free products of \textit{finite} groups. In this case we have
$\Omega=\Omega_\infty$. Throughout the whole subsection we do not need the
assumption that the $\mu_i$'s are symmetric. For any
$\omega_1=x_1x_2\dots,\omega_2=y_1y_2\dots\in\Omega_\infty$ with $\omega_1\neq
\omega_2$, we define the \textit{confluent} $\omega_1\wedge \omega_2$ of
$\omega_1$ and $\omega_2$ as the word $x_1\dots x_k$ of maximal length with
$x_i=y_i$ for all $1\leq i\leq k$. If $x_1\neq y_1$, then $\omega_1 \wedge \omega_2:=e$. The metric on the boundary $\Omega_\infty$ is defined by
$$
d^{\mathrm{fin}}_{\Omega}(\omega_1,\omega_2) := \alpha^{\Vert \omega_1\wedge\omega_2\Vert}
$$
for any arbitrary but fixed $\alpha\in (0,1)$. With respect to this metric on
$\Omega_\infty$ we can define analogously to (\ref{def:BD}) and (\ref{def:HD}) the upper box-counting dimension
$\overline{\mathrm{BD}^{\mathrm{fin}}}(\Omega')$, the box-counting dimension
$\mathrm{BD}^{\mathrm{fin}}(\Omega')$ and the Hausdorff dimension
$\mathrm{HD}^{\mathrm{fin}}(\Omega')$ for any $\Omega'\subseteq
\Omega_\infty$. We set $\mathcal{F}_i^+(\lambda):=\mathcal{F}^+_i(\lambda|1)$ and
define the matrix $M=\bigl(m(i,j)\bigr)_{i,j\in\mathcal{I}}$ by
$$
m(i,j):=
\begin{cases}
\mathcal{F}_j^+(\lambda), & \textrm{if } i\neq j,\\
0,&\textrm{if } i=j.
\end{cases}
$$
Since $M$ is irreducible and has non-negative entries, the Perron--Frobenius
eigenvalue exists and is denoted by $\theta$.
\par
Furthermore, define the matrix $D=\bigl(d(i,j)\bigr)_{i,j\in\calI}$ by $d(i,j):=|\Gamma_j|-1$, if $i\neq j$, and $d_{i,i}:=0$, and denote by $\varrho$ its Perron--Frobenius eigenvalue. With this notation we get:
\begin{Cor}\label{Cor:finite-groups}
$$
\mathrm{BD}^{\mathrm{fin}}(\Lambda) = \mathrm{HD}^{\mathrm{fin}}(\Lambda) = -\frac{\log\theta}{\log\alpha}
\quad \textrm{ and } \quad
\mathrm{BD}^{\mathrm{fin}}(\Omega) = \mathrm{HD}^{\mathrm{fin}}(\Omega) = -\frac{\log\varrho}{\log\alpha}.
$$
\end{Cor}
\begin{flushright}
$\Box$
\end{flushright}
Let us remark that, in the case of $\Gamma=\Gamma_1\ast\Gamma_2$ with $|\Gamma_1|=|\Gamma_2|<\infty$, we get the following explicit formulas for the dimensions:
\begin{eqnarray*}
\mathrm{BD}^{\mathrm{fin}}(\Lambda) & = & \mathrm{HD}^{\mathrm{fin}}(\Lambda) = -\frac{\log \sqrt{\mathcal{F}_1^+(\lambda)\mathcal{F}^+_2(\lambda)}}{\log\alpha}\ \textrm{ and}\\
\mathrm{BD}^{\mathrm{fin}}(\Omega) & = & \mathrm{HD}^{\mathrm{fin}}(\Omega) = -\frac{\log \sqrt{\bigl(|\Gamma_1|-1 \bigr)\bigl(|\Gamma_2|-1 \bigr)}}{\log\alpha}.
\end{eqnarray*}
\par
\addtocounter{Th}{1}
\begin{Example}
Consider $\Gamma=(\Z/3\Z)\ast(\Z/2\Z)$, where $\Z/3\Z=\{e_1,a,a^2\}$ and $\Z/2\Z=\{e_2,b\}$. We choose $\mu(a)=p\in(0.1,0.7)$, $\mu(a^2)=q\in (0,0.9-p)$ and $\mu(b)=1-p-q$. We set $\alpha:=1/2$ and $\lambda=1.005$. Let us note that this choice of the parameters $p$ and $q$ lead to $R\geq 1.005$, which can be verified by numerical evaluation.
For instance, in \cite[Section 3.6.1]{gilch} the required generating functions are computed. In Figure  \ref{fig:boxdim-z2-z3} we can see the behaviour of $\mathrm{HD}^{\mathrm{fin}}(\Lambda)$ with $\lambda=1.005$ in dependence of the parameters $p$ and $q$. The Hausdorff dimension of the whole space of ends is $ 0.5$; compare with Example \arabic{helpcounter1}.\arabic{helpcounter2}.
\begin{figure}
\begin{center}
\includegraphics[width=6cm]{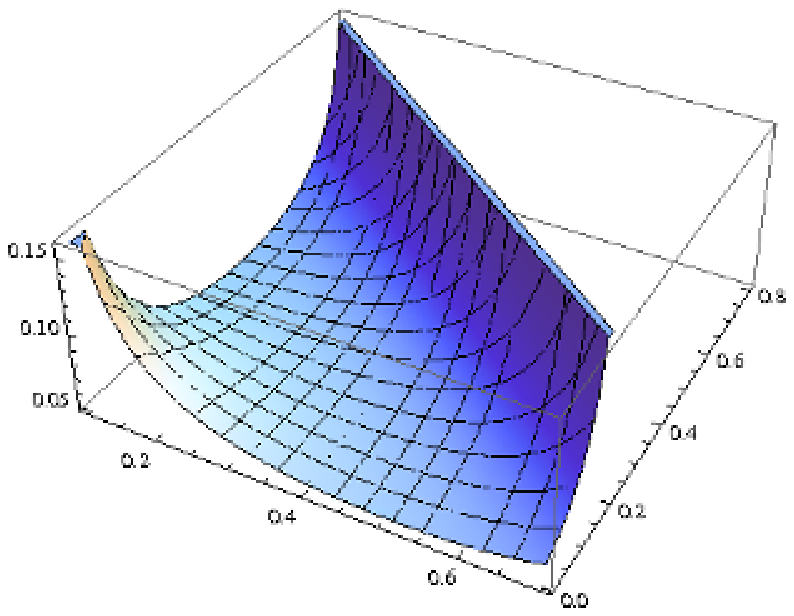}
\end{center}
\caption{Hausdorff dimension $\mathrm{HD}(\Lambda)$ of the branching random
walk on \mbox{$(\Z/3\Z)\ast(\Z/2\Z)$} with $\lambda=1.005$ in dependence of $p$ and $q$.}
\label{fig:boxdim-z2-z3}
\end{figure}
\end{Example}

\subsection{Free Products by Amalgamation of Finite Groups}
An important generalization of free products are free products by amalgamation (of finite groups). Let $\Gamma_1,\dots,\Gamma_r, H$, be finite groups such that each group $\Gamma_{i}$ contains  a subgroup $H_{i}$ that is isomorphic to $H$. Let $\phi_{i}: H_{i}\to H$ be  an isomorphism for each $i\in\{1,\ldots,r\}$. Moreover, let $S_{i}$ be a generating set of $\Gamma_{i}$ and $R_{i}$ its  relations. 
The \textit{free product by amalgamation} with respect to the subgroup $H$ is defined by
\begin{eqnarray*}
\Gamma_H &:=& \Gamma_1 \ast_H \Gamma_2 \ast_H \dots \ast_H \Gamma_r\\
&:=& \langle S_{1},\ldots, S_{r}\mid R_{1},\ldots, R_{n}, \phi_{j}^{-1}(\phi_{i}(a))=a~\forall a\in H_{i}~\forall  i,j\in\mathcal{I}\rangle.
\end{eqnarray*}
For $i\in\calI$, the quotient $\Gamma_i/H_i$ consists of all left co-sets of the form \mbox{$x_iH_i = \{x_ih \mid h \in H_i\}$,} where $x_i\in\Gamma_i$. We fix a set of representatives $\mathcal{R}_i:=\{g_{i,1} = e_i, g_{i,2},\dots, g_{i,n_i}\}$ for the elements of $\Gamma_i /H_i$, that is, for each $y_i\in\Gamma_i$ there is a unique $g_{i,k}\in \mathcal{R}_i$ with $y_i\in g_{i,k}H_i$. 
We write  $\hat\tau(x)=i$ if $x\in \mathcal{R}_i\setminus \{e_i\}$.
The amalgam $\Gamma_H $ consists of all finite words of the form
\begin{equation}\label{def:amalgam-words}
x_1 x_2 \dots x_nh
\end{equation}
with $n\in \N_0$, $x_i\in \bigcup_{j\in\calI} \mathcal{R}_j\setminus \{e_j\}$
and $h\in H$ such that $\hat\tau(x_i)\neq \hat\tau(x_{i+1})$. Here w.l.o.g. we
may identify $h$ with $\phi_1^{-1}(h)$, and $e$ denotes again the empty word.
Let $\Omega$ be the set of all ends of $\Gamma_H$, which consists of all infinite words of the form $w_1w_2\ldots\in \Bigl(\bigcup_{i\in\calI} \mathcal{R}_i\setminus \{e_i\}\Bigr)^{\N}$ such that $\hat\tau(w_i)\neq \hat\tau(w_{i+1})$ for all $i\in\N$. For any $\omega_1=x_1x_2\dots,\omega_2=y_1y_2\ldots\in\Omega$ with $\omega_1\neq \omega_2$, we define again the \textit{confluent} $\omega_1\wedge \omega_2$ of $\omega_1$ and $\omega_2$ as the word $x_1\dots x_k$ of maximal length with $x_i=y_i$ for all $1\leq i\leq k$. If $x_1\neq y_1$, then $\omega_1\wedge \omega_2:=e$. Again we can define a metric on the boundary $\Omega$:
$$
d^{(H)}_{\Omega}(\omega_1,\omega_2) := \alpha^{\Vert \omega_1\wedge \omega_2\Vert}
$$
for any $\alpha\in (0,1)$. With respect to this metric on $\Omega$ we can define analogously to (\ref{def:BD}) and (\ref{def:HD}) the upper box-counting dimension $\overline{\mathrm{BD}^{(H)}}(\Omega')$, the box-counting dimension $\mathrm{BD}^{(H)}(\Omega')$ and Hausdorff dimension $\mathrm{HD}^{(H)}(\Omega')$ for any $\Omega'\subseteq \Omega$.
\par
Suppose we are given symmetric probability measures $\mu_i$ on the groups
$\Gamma_i$ and numbers $\alpha_{i}>0$ such that $\sum_{i\in\mathcal{I}}
\alpha_{i}=1$. The random walk on $\Gamma_H$ is then governed by
$$
\mu(x) :=
\begin{cases}
\alpha_i \mu_i(x), & \textrm{if } x\in\Gamma_i\setminus H_i,\\
\sum_{i\in\calI} \alpha_i \mu_i\bigl( \phi^{-1}_i(\phi_1(x)) \bigr), & \textrm{if } x\in H_1,\\
0, & \textrm{otherwise}.
\end{cases}
$$
For $g_i\in \mathcal{R}_i$, denote by $T_{g_iH}$ the stopping time of the first visit of the set $g_iH_i$.
We introduce the following generating functions:
$$
F_H(gh|z) := \sum_{n\geq 0} \Prob\bigl[T_{gH}=n, X_n=gh \mid X_0=e \bigr]\,z^n,
$$
where $g\in \bigcup_{i\in\calI} \mathcal{R}_i\setminus \{e_i\}$, $h\in H_i$ and
$z\in\mathbb{C}$. By symmetry of the $\mu_i$'s, we have $F_H(gh|z)\leq F(e,gh|z)<1$; compare with (\ref{equ:F<1}).
Conditioning on the first step of the random walk, we get
\begin{equation}\label{equ:F_H}
\begin{array}{rcl}
F_H(gh|z) & = & \mu(gh) z + \sum_{g_0\in \Gamma_{\tau(g)}\setminus gH_{\tau(g)}} \mu(g_0)z F_H(g_0^{-1}gh|z) +  \\[2ex]
&& \quad 
\sum_{i\in\calI\setminus\{\tau(g)\}} \sum_{g_0\in\Gamma_i}\mu(g_0) z \sum_{h_0\in H_i} F_H(g_0^{-1}h_0|z) F_H(h_0^{-1}gh|z).
\end{array}
\end{equation}
Since there are only \textit{finitely} many functions $F_H(\cdot|z)$, one can compute $F_H(\cdot|z)$ by solving the finite system of quadratic equations (\ref{equ:F_H}). We define also
$$
\mathcal{F}^{(H)}_i(z) := \sum_{\substack{g\in\mathcal{R}_i\setminus\{e_i\},\\
    h\in H_i}} F_H(gh|z)
$$
and the matrix  $N=\bigl(n(i,j)\bigr)_{i,j\in\calI}$ with entries
$$
n(i,j):=
\begin{cases}
\mathcal{F}^{(H)}_j(\lambda), & \textrm{if } i\neq j,\\
0,&\textrm{if } i=j.
\end{cases} 
$$
We denote by $\theta_H$ the Perron--Frobenius eigenvalue of $N$. Furthermore, we denote by $\varrho_H$ the Perron--Frobenius eigenvalue of the matrix $D_H=\bigl(d_H(i,j)\bigr)_{i,j\in\mathcal{I}}$, which is defined by
$$
d_H(i,j):=
\begin{cases}
[\Gamma_j: H_j]-1,& \textrm{if } i\neq j,\\
0 & \textrm{if } i=j.
\end{cases}
$$
Finally, we can state the following formulas for the dimensions:
\begin{Cor}\label{Cor:amalgam}
$$
\mathrm{BD}^{(H)}(\Lambda) = \mathrm{HD}^{(H)}(\Lambda) = -\frac{\log\theta_H}{\log\alpha}
\quad \textrm{ and } \quad 
\mathrm{BD}^{(H)}(\Omega) = \mathrm{HD}^{(H)}(\Omega) = -\frac{\log\varrho_H}{\log\alpha}.
$$
\end{Cor}
\par
\addtocounter{Th}{1}
\begin{Example}
Consider the amalgam $(\Z/6\Z) \ast_{\Z/2\Z}(\Z/6\Z)$. Hence, let
\mbox{$\Gamma_{1}=\langle a \mid a^{6}=e_1\rangle$,} $\Gamma_{2}=\langle b \mid
b^{6}=e_2\rangle$, and $H=\langle c \mid c^{2}=e_H\rangle$, where $e_H$ is the
identity in $H$. The isomorphisms are defined through $\phi_{1}(a^{3})=c=\phi_{2}(b^{3})$. Eventually,
$$
(\Z/6\Z) \ast_{\Z/2\Z}(\Z/6\Z)= \langle a,b \mid a^{6}=b^{6}=e,
a^{3}=b^{3}\rangle.
$$
We set $\mu_1(a)=\mu_1(a^5)=\mu_2(b)=\mu_2(b^5)=1/2$, $\alpha_1=\alpha_2=1/2$
and consider the distance with base $\alpha=1/2$. The system (\ref{equ:F_H}) becomes then
\begin{eqnarray*}
F_H(a|z) & = & \frac{z}{4} + \frac{z}{4}F_H(a^2|z) + \frac{z}{2}\bigl( F_H(a|z)^2 + F_H(a^2|z)^2\bigr),\\
F_H(a^2|z) & = & \frac{z}{4}F_H(a|z) + \frac{z}{2}\bigl( F_H(a|z) F_H(a^2|z) + F_H(a^2|z) F_H(a|z)\bigr).
\end{eqnarray*}
Observe that $F_H(a|z)=F_H(a^5|z)$ and $F_H(a^2|z)=F_H(a^4|z)$. The Hausdorff dimension of the branching random walk is then given by
$$
\mathrm{HD}^{(H)}(\Lambda) = \frac{\log \bigl(2F_H(a|\lambda) + 2F_H(a^2|\lambda) \bigr)}{\log 2},
$$
while $\mathrm{HD}^{(H)}(\Omega)=1$. The behaviour of $\mathrm{HD}^{(H)}(\Lambda)$ in function of $\lambda$ is qualitatively the same as in Figure \ref{fig:hd-z2-z3}.
\end{Example}

\section{Proofs}\label{sec:proofs}

\subsection{Proof of Theorem \ref{thm:phase-transition}}\label{sec:proof1}
We first introduce some preliminary results on BRW. 
Using the description of a tree-indexed random walk it is easy to see  that the distribution of the location of some particle in
generation $n$ has the same distribution as the location of a (non-branching)
random walk on $\Gamma$ after $n$ steps, see \cite{benjamini:94}.
\begin{Lemma}\label{lemma:n-step}
Let $v\in \T$ with $|v|=n$ for some $n\geq 1$. Then,
$$\P[S_{v}=y]=P[X_{n}=y]=\mu^{(n)}(y).$$
\end{Lemma}
The following lemma will be used several times in our proofs. It gives a
formula for the expected number of elements frozen in a set $M$, in  the coloured branching random walk. This observation can be found for example in \cite{MV} or \cite[Lemma 1]{lalley2000}. Nevertheless, we give a short proof since it is one of the essential points where the generating function $F(\cdot,\cdot|z)$ intervenes.
\begin{Lemma}\label{lemma:Z-infty-formel}
For any $M\subseteq \Gamma$, we have $\mathbb{E}\bigl[Z_\infty(M)\bigr]=F(e,M|\lambda)$. 
\end{Lemma}
\begin{proof}
For any $v\in \T$, let  $\langle v_0=\ro, v_1, \ldots, v_{|v|}=v\rangle$ be the unique geodesic from $\ro $ to $v$. Now, we define for any $n\in\N$ 
$$
\mathrm{Fr}_{v}^{(n)}:=
\left\{
\begin{array}{ll}
  1, & \mbox{if } v\in M \textrm{ and } v_{i}\notin M~\forall i\leq n-1,  \\
  0, & \mbox{otherwise.}     
\end{array}\right.
$$
 In words, $\mathrm{Fr}^{(n)}_v$ is the number of particles being frozen in $v$ at time $n$.   Using the well-known fact that $E[|\T_{n}|]=\lambda^{n}$ we obtain
\begin{eqnarray*}
\mathbb{E}\Bigl[ \sum_{v\in\T_{n}}\mathrm{Fr}^{(n)}_v \Bigr] &= &\sum_{k\geq 1} \mathbb{E}\Bigl[ \sum_{v\in\T_{n}}\mathrm{Fr}^{(n)}_{v} \Bigl| |\T_{n}| =k \Bigr] \P[|\T_{n}| =k] \\
& = & \sum_{k\geq 1} \Prob\bigl[X_n\in M,\forall~  m\leq n-1: X_m\notin M\bigr]\,k\,\P[|\T_{n}| =k]\\
& =& \Prob\bigl[X_n\in M,\forall~  m\leq n-1: X_m\notin M\bigr]\,\lambda^n.
\end{eqnarray*}
Summing over $n$ finishes the proof.
\end{proof}

The proof of Theorem \ref{thm:phase-transition} splits up into the proofs of the following Propositions \ref{prop:xi-bigger-1}, \ref{prop:lambda-bigger-1-infty} and \ref{prop:xi-equal-1}.

Recall from the definition of $\Omega_i^{(0)}$ and $\Omega_i$ that
$\Omega_i^{(0)}\subseteq \Omega_i\subseteq \Omega$.
\begin{Prop}\label{prop:xi-bigger-1}
Ends of $\Omega_i^{(0)}$ occur in $\Lambda$ with positive probability if and only if $\xi_i(\lambda)>1$, that is, $\Prob\bigl[\Lambda \cap \Omega_i^{(0)}\neq \emptyset\bigr]>0$
if and only if $\xi_i(\lambda)>1$.
\end{Prop}
\begin{proof}
It is convenient to work with the  \textit{coloured} branching random walk.  In
fact, the idea of the proof is to define an embedded Galton--Watson process
that counts the number of particles that hit $\Gamma_i$, where $\xi_i(\lambda)$
will be the growth parameter.
%

We start the BRW with one particle in $e=e_{i}$.
The first generation of the branching process is formed by those particles that
are frozen in  $\Gamma_i^\times$. Let us check that the number of those
particles is almost surely finite. Since $\mu$ has finite
support every particle visiting $\Gamma_i^\times$ has to pass through
$\mathrm{supp}(\mu_i)$. Hence,
$Z_{\infty}(\Gamma_i^\times)=Z_{\infty}(\mathrm{supp}(\mu_i))$, which is almost
surely finite since the BRW is transient. 
The second generation of the branching process is constructed as follows. For each particle frozen in some $x\in \Gamma_i^\times$ we start a new BRW where each particle when reaching $\Gamma_i \setminus \{x\}$ is frozen. Now, the second generation of the branching process consists of all these new frozen particles. Further generations are constructed inductively in the same way. Let $\psi_{n}$ be the number of particles of this process at generation $n$. Obviously, $(\psi_{n})_{n\geq 0}$ turns out to be a Galton--Watson process with mean
$$
m_i=\mathbb{E}\bigl[Z_\infty\bigl(\mathrm{supp}(\mu_i)\bigr)\bigr]=F\bigl(e,\mathrm{supp}(\mu_i)|\lambda\bigr)=\xi_i(\lambda).
$$
Hence,  this Galton--Watson process survives with positive probability if and only if
$\xi_i(\lambda)>1$; see e.g. \cite[Theorem 6.1]{harris63}. As a consequence, we
have that $\Gamma_{i}$ is visited infinitely many times with positive
probability if $\xi_i(\lambda)>1$. That is,  $\Prob\bigl[\Lambda \cap
\Omega_i^{(0)}\neq \emptyset\bigr]>0$ if $\xi_i(\lambda)>1$. On the other hand,
$\xi_i(\lambda)\leq 1$ implies that $\Gamma_{i}$ is almost surely visited only
for a finite number of times and hence $\Prob\bigl[\Lambda \cap \Omega_i^{(0)}\neq \emptyset\bigr]=0.$
\end{proof}

The next step is to  show $|\Lambda \cap \Omega_i|=\infty$ if $\xi_{i}(\lambda)>1$.
\begin{Prop}\label{prop:lambda-bigger-1-infty}
If $\xi_i(\lambda)>1$ then there are almost surely infinitely many cosets 
$x\Gamma_i$, where the branching random walk accumulates. That is, the set 
$$
\bigl\lbrace
x\in\Gamma \,\bigl|\, \tau(x)\neq i, x\Omega_i^{(0)} \cap \Lambda\neq
\emptyset\bigr\rbrace
$$
is almost surely infinite. 
\end{Prop}
\begin{proof}
We construct the family tree $\T$ of the BRW with branching distribution $\nu$ in the following way. We start with one geodesic line  $v_{\infty}=\langle \ro, v_{1}, v_{2}, \ldots\rangle$ and attach to each of the vertices independent copies of Galton--Watson trees where the distribution of the first generation is $\tilde \nu(k)=\nu(k+1)$ for $k\geq 0$ and $\nu$ for the other generations. 
The trajectory along $v_\infty$ has the same distribution as a non-branching random walk, compare with Lemma \ref{lemma:n-step}. Hence,  $S_{v_{n}}$ converges almost surely to a random infinite word
$g_\infty=g_1g_2\ldots\in\Omega_\infty$  as $n\to\infty$; here we mean convergence in the sense that the
block length of the common prefix of the location of $S_{v_n}$ and $g_\infty$ tends to
infinity. Moreover, we define the
random indices $n_1:=\min\{m\in\N \mid g_m\in\Gamma_i\}$, and recursively $n_k:=
\min\{m\in\N \mid m>n_{k-1}, g_m\in\Gamma_i\}$. Note that these indices are
almost surely finite; see e.g. \cite[Section 7.I]{gilch:07}. Denote by $\hat
v_k$ the first vertex in $v_\infty$ with $\hat v_k=g_1\dots g_{n_k}$.
Let $B_{k}$ be the set of offspring of  $\hat v_{k}=v_s$ different from
$v_{s+1}$ and denote by $\Lambda_{v}$ the set of accumulation points of the
descendants of some $v\in\T$. Moreover, we define $A_{k}$ as the event that
$\Lambda_{v}\cap S_v\Omega_{i}^{(0)}\neq\emptyset$ for some $v\in B_{k}$ with
$\tau(v)=i$. Observe that the events $A_{k}$ are i.i.d. since transitivity yields
$\Prob[\Lambda_{v}\cap S_v\Omega_{i}^{(0)}\neq\emptyset]=\Prob[\Lambda \cap
\Omega_i^{(0)}\neq \emptyset]$ for every $v\in\T$. Now, due to Proposition  \ref{prop:xi-bigger-1}
and the fact that 
\begin{eqnarray*}
P[B_{k}\neq \emptyset, \exists v\in B_k: \tau(S_v)=i] &= &
\bigl(1-\nu(1)\bigr)\cdot \Prob[v\in B_k:\tau(S_v)=i\mid B_{k}\neq \emptyset]\\
&\geq &\bigl(1-\nu(1)\bigr)\cdot\alpha_i
>0
\end{eqnarray*}
we have $P[A_{k}]\geq c$ for all $k$ and some $c>0$. Eventually, the Lemma of  Borel--Cantelli yields that an infinite number of $A_{k}$'s occurs almost surely. 
\end{proof}

It remains to treat the critical and subcritical case $\xi_i(\lambda)\leq 1$ in order to complete  the proof of Theorem \ref{thm:phase-transition}.
\begin{Prop}\label{prop:xi-equal-1}
If $\xi_i(\lambda)\leq 1$ then $\Prob[\Lambda \cap \Omega_i\neq \emptyset]=0$.
\end{Prop}
\begin{proof}
Due to Proposition \ref{prop:xi-bigger-1} we have that 
$\Prob\bigl[ \Lambda \cap x\Omega_i^{(0)}\neq \emptyset\bigr]=0$ for all $x\in
\Gamma$: indeed, each $x\in\Gamma$ is almost surely visited finitely often;
each particle, which hits $x$, starts its own branching random walk at $x$ and
each of these branching random walk hits $x\Omega_i^{(0)}$ only finitely often
with probability one. Since $$\Lambda\cap \Omega_{i}=\biguplus_{x\in\Gamma: \tau(x)\neq i}(\Lambda\cap x \Omega_{i}^{(0)})$$ we conclude
$$ \Prob\bigl[ \Lambda \cap \Omega_i\neq \emptyset\bigr]=\sum_{x\in\Gamma: \tau(x)\neq i} \Prob\bigl[ \Lambda \cap x\Omega_i^{(0)}\neq \emptyset\bigr]=0.$$
\end{proof}

\subsection{Proof of Theorem \ref{thm:boxdim} and Corollary \ref{cor:omega-i-do-not-contribute}}

First, we show  that the proposed
formula for the dimension is an upper bound for the upper box-counting
dimension; see Proposition \ref{prop:upperbound} in Subsection
\ref{subsec:upperbound}. In the second step we show that the proposed formula
is also a lower bound for the Hausdorff dimension of $\Lambda$; see Corollary
\ref{cor:boxdim-lambda} in Subsection \ref{subsec:lower-bound}. Finally, this will imply  the proof of Theorem \ref{thm:boxdim} and Corollary \ref{cor:omega-i-do-not-contribute}.

%

\subsubsection{Upper Bound for the Box-Counting Dimension}

\label{subsec:upperbound}
In this part we show that $\log z^\ast/\log\alpha$ is an upper bound for
$\overline{\mathrm{BD}}(\Lambda)$. 
To this end we introduce the following  notation: for $n\in\N$, we denote by 
$$
\mathcal{H}_n :=\bigl\lbrace x\in \Gamma \,\bigl|\, l(x)=n, x \textrm{ is
 visited by the branching random walk} \bigr\rbrace
$$
the set of visited sites at graph distance $n$. An important observation is that 
for each end $\omega\in\Lambda$ and every $m\in\N$, the branching random walk
has to visit at
least one vertex $x_\omega\in\mathcal{H}_m$, where  $x_\omega$ is in the
$\omega$-component of $\mathcal{X}\setminus B_{m-1}$. Thus,
$$
\Lambda \subseteq \bigcup_{x\in\mathcal{H}_m} \bigl\lbrace \omega\in \Omega
\mid x \textrm{ lies in the } \omega\textrm{-component of } \mathcal{X}\setminus
B_{m-1}\bigr\rbrace.
$$ 
This implies that $\Lambda$ can be covered by $|\mathcal{H}_m|$ balls of radius
$\alpha^m$. Our strategy for the upper bound is to study the limit behaviour of
$\mathbb{E} |\mathcal{H}_m|^{1/m}$ first and then the resulting  limit behaviour of  $|\mathcal{H}_m|^{1/m}$ as $m\to\infty$; see Lemma \ref{lemma:upperbound}. This
will eventually  lead to the proposed upper bound for $\overline{\mathrm{BD}}(\Lambda)$; see Proposition \ref{prop:upperbound}.
\par
Observe that $x\in\mathcal{H}_m$ if and only if $Z_\infty(x)\geq 1$. Therefore, by Lemma \ref{lemma:Z-infty-formel},
$$
1\leq \mathbb{E}|\mathcal{H}_m|\leq \sum_{x\in\Gamma: l(x)=m} \mathbb{E}Z_\infty(x) =
\sum_{x\in\Gamma: l(x)=m} F(e,x|\lambda)=:H_m.
$$
We have that $H_{m+n}\leq H_{m} H_{n}$ and hence Fekete's lemma implies that $\lim_{n\to \infty} H_{m}^{1/m}$ exists. 
Recall the definitions of $\mathcal{F}(\lambda|z)=\sum_{m\geq 0} H_m\,z^m$, $\mathcal{F}_i^+(\lambda|z)$ and $\mathcal{F}_i(\lambda|z)$ in (\ref{equ:scriptF}), (\ref{equ:scriptFi+}) and (\ref{equ:scriptFi}).  Due to (\ref{equ:scriptF2}) we get the equation
$$
\mathcal{F}_i(\lambda|z)= \mathcal{F}_i^+(\lambda|z) \bigr(\mathcal{F}(\lambda|z)-\mathcal{F}_i(\lambda|z)\bigr),
$$
or equivalently
$$
\mathcal{F}_i(\lambda|z)=\mathcal{F}(\lambda|z) \frac{\mathcal{F}_i^+(\lambda|z)}{1+\mathcal{F}_i^+(\lambda|z)}.
$$
Hence,
$$
\mathcal{F}(\lambda|z)=1+\sum_{i\in\mathcal{I}}\mathcal{F}_i(\lambda|z)= 1+ \mathcal{F}(\lambda|z)\sum_{i\in\mathcal{I}} \frac{\mathcal{F}_i^+(\lambda|z)}{1+\mathcal{F}_i^+(\lambda|z)},
$$
or equivalently
\begin{equation}\label{equ:F-equation}
\mathcal{F}(\lambda|z)=\frac{1}{1-\sum_{i\in\mathcal{I}} \frac{\mathcal{F}_i^+(\lambda|z)}{1+\mathcal{F}_i^+(\lambda|z)}}.
\end{equation}
This equation holds for every $z\in\mathbb{C}$ with $|z|<R(\mathcal{F})$, where
$R(\mathcal{F})$ is the radius of convergence of $\mathcal{F}(\lambda|z)$. 
Since $$1\leq \lim_{m\to\infty} H_m^{1/m} = 1/R(\mathcal{F}),$$ we have
\begin{equation}\label{equ:Fless1}
R(\mathcal{F})\leq 1.
\end{equation}
 In order to determine $R(\mathcal{F})$ 
we have to find -- by Pringsheim's Theorem -- the smallest
singularity point on the positive $x$-axis of $\mathcal{F}(\lambda|z)$. This
smallest singularity point is either one of  the radii of convergence
$R(\mathcal{F}_i^+)$ of the
functions $\mathcal{F}_i^+(\lambda|z)$ or the smallest real positive number
$z^\ast$ with 
\begin{equation}
\sum_{i\in\mathcal{I}}
\frac{\mathcal{F}_i^+(\lambda|z^\ast)}{1+\mathcal{F}_i^+(\lambda|z^\ast)}=1.
\end{equation}
The next two lemmas imply that in fact $R(\mathcal{F})=z^\ast$.
\begin{Lemma}\label{prop:conv-radius}
$R(\mathcal{F}) \in (0,1)$
\end{Lemma}
\begin{proof}\footnote{This short proof was suggested by the referee.}
The fact that $R(\mathcal{F})>0$ follows from the fact that the Cayley graph grows not faster than exponentially. To see that $R(\mathcal{F})<1$ recall that Equation (\ref{gl-equations}) states that the generating functions $F(e,x|z)$ and $G(e,x|z)$ are comparable, i.e., $G(e,x| \lambda)=F(e,x|\lambda)G(e,e|\lambda)$. Hence, for some $C>0$ we have for all $m\in\N$ that
$$\sum_{x: l(x)\leq m} F(e,x|\lambda)\geq C \sum_{x: l(x)\leq m} G(e,x|\lambda).$$
The sum on the right hand side is the expected number of visits of the BRW in the
ball $B_{m}$, the set of vertices $x\in\Gamma$ with $l(x)\leq m$. As we assumed
the random walk to be of nearest neighbour type  all particles up to generation
$m$ must be contained in the ball $B_{m}$. The expected population size at time
$m$ is just $\lambda^{m}$ which eventually  implies that  $H_{m}$ grows
exponentially fast, since $\lim_{m\to\infty} H_m^{1/m}$ exists and is at least $1$.

\end{proof}
\begin{Lemma}\label{lemma:conv-radius2}
For all $i\in\calI$, $R(\mathcal{F})=z^\ast < R(\mathcal{F}_i^+)$.
\end{Lemma}
\begin{proof}
Let us first consider the case  $\xi_i(\lambda)<1$, where we obtain
\begin{equation}\label{equ:fi+}
\begin{array}{rcl}
\mathcal{F}_i^+(\lambda|1) &=& \sum_{x\in\Gamma_i^\times} F_i\bigl(e_i,x|\xi_i(\lambda)\bigr)\\[2ex]
& = & \frac{1}{G_i\bigl(e_i,e_i\bigl|\xi_i(\lambda)\bigr)} \sum_{x\in\Gamma_i} G_i\bigl(e_i,x\bigl|\xi_i(\lambda)\bigr) -1 \\[2ex]
&=&
\frac{1}{G_i\bigl(e_i,e_i|\xi_i(\lambda)\bigr) \bigl(1-\xi_i(\lambda)\bigr)} -1<\infty.
\end{array}
\end{equation}
Hence,  $\xi_i(\lambda)<1$ implies $R(\mathcal{F}_i^+)\geq 1>R(\mathcal{F})$. 
In the case of $\xi_i(\lambda)\geq 1$ the claim follows from the following inequality:
\begin{equation}\label{equ:limsup-inequ}
\frac{1}{R(\mathcal{F})}=\limsup_{n\to\infty} \biggl(\sum_{\substack{x\in\Gamma:\\ l(x)=n}} F(e,x|\lambda)\biggr)^{1/n}
> 
\limsup_{n\to\infty} \biggl(\sum_{\substack{x\in\Gamma_1:\\ l(x)=n}} F(e,x|\lambda)\biggr)^{1/n}=\frac{1}{R(\mathcal{F}_1^+)}.
\end{equation}
In order to prove (\ref{equ:limsup-inequ}) we define for 
 $n\in\N$
$$
a_n:=\log  \sum_{\substack{x\in\Gamma_1:\\ l(x)=n}} F(e,x|\lambda).
$$
We have that $a_n\geq 0$ since
$$
 \sum_{\substack{x\in\Gamma_1:\\ l(x)=n}} F(e,x|\lambda)= \sum_{\substack{x\in\Gamma_1:\\ l(x)=n}} F_1\bigl(e_1,x|\xi_1(\lambda)\bigr)
\geq \sum_{\substack{x\in\Gamma_1:\\ l(x)=n}} F_1\bigl(e_1,x|1)
\geq \Prob[T_{S_1(n)}<\infty]=1,
$$
where $S_1(n):=\{x\in\Gamma_1\mid l(x)=n\}$ and $T_M$ is the stopping time for the random walk on $\Gamma_1$ (governed by $\mu_1$) of the first visit of a set $M\subseteq \Gamma_1$. Furthermore, $(a_n)_{n\in\N}$ is a subadditive sequence, that is, $a_m+a_n\geq a_{m+n}$ for all $m,n\in\N$.
By  Fekete's Lemma, the limit $\lim_{n\to\infty} a_{n}/n$ exists and is equal to $\inf_{n\in\N} a_{n}/n$, hence
$$
\lim_{n\to \infty} \biggl(\sum_{\substack{x\in\Gamma_1:\\ l(x)=n}} F(e,x|\lambda)\biggr)^{1/n}= \frac{1}{R(\mathcal{F}_1^+)}
=\inf_{n\in\N} \biggl(\sum_{\substack{x\in\Gamma_1:\\ l(x)=n}} F(e,x|\lambda) \biggr)^{1/n}.
$$
The last equation implies that
$$
\biggl(\sum_{\substack{x\in\Gamma_1:\\ l(x)=n}} F(e,x|\lambda)\biggr)^{1/n}\geq \frac{1}{R(\mathcal{F}_1^+)}=:q_1 \quad \forall n\in\N.
$$ 
Observe that $\sum_{x\in\Gamma_2:l(x)=1}F(e,x|\lambda)\geq \xi_2(\lambda)$.
Then, for all $n\in\N$: 
\begin{eqnarray*}
H_n=\sum_{\substack{x\in\Gamma:\\ l(x)=n}} F(e,x|\lambda) & = & \sum_{k=1}^n \sum_{\substack{x=x_1\dots x_k\in\Gamma:\\ l(x)=n}} \prod_{j=1}^k F(e,x_j|\lambda)\\
&\geq & \sum_{k=1}^{\lfloor n/2\rfloor} \sum_{\substack{x_1,\dots, x_k\in\Gamma_1:\\  l(x_1)+\dots l(x_k) +k=n}} \xi_2(\lambda)^k \prod_{j=1}^k F(e,x_j|\lambda)\\
&\geq & \sum_{k=1}^{\lfloor n/2\rfloor} \sum_{\substack{n_1,\dots, n_k\in\N:\\  n_1+\dots +n_k+k=n}} q_1^{n_1} \xi_2(\lambda) q_1^{n_2} \xi_2(\lambda) q_1^{n_3} \cdots \xi_2(\lambda) q_1^{n_k} \xi_2(\lambda)\\
&\geq & \sum_{k=1}^{\lfloor n/2\rfloor} q_1^{n-k} \xi_2(\lambda)^{k} \binom{n-2k+k-1}{k-1}.
\end{eqnarray*}
In the last inequality the binomial coefficients arise as follows: we think of counting the number of possibilities of placing $n-k$ (undistinguishable) balls into $k$ urns, where each urn should at least contain one ball.
We note that $n-k-2\geq \lfloor n/2\rfloor-1$ for all $k\leq \lfloor n/2\rfloor-1$. Therefore, with the help of the Binomial theorem we obtain:
\begin{eqnarray*}
H_n & \geq &  q_1^n  \sum_{k=0}^{\lfloor n/2\rfloor-1} \Bigl(\frac{\xi_2(\lambda)}{q_1}\Bigr)^{k+1}  \binom{n-k-2}{k}\\
&\geq &
q_1^{n-1} \xi_2(\lambda) \sum_{k=0}^{\lfloor n/2\rfloor-1} \Bigl(\frac{\xi_2(\lambda)}{q_1}\Bigr)^k \binom{\lfloor n/2\rfloor -1}{k}
\geq  q_1^{n-1}\xi_2(\lambda) \Bigl(1+\frac{\xi_2(\lambda)}{q_1}\Bigr)^{\lfloor n/2\rfloor-1}.
\end{eqnarray*}
Taking $n$-th roots on both sides  and letting $n\to\infty$ yields 
\begin{equation}\label{equ:rF-rFi}
\liminf_{n\to \infty}\Bigl(\sum_{\substack{x\in\Gamma:\\ l(x)=n}} F(e,x|\lambda)\Bigr)^{1/n} 
\geq \frac{1}{R(\mathcal{F}_1^+)} \sqrt{1+\frac{\xi_2(\lambda)}{q_1}}>\frac{1}{R(\mathcal{F}_1^+)}.
\end{equation}
\end{proof}

The next lemma gives an almost sure upper bound for $|\mathcal{H}_m|^{1/m}$ as $m\to\infty$. Its proof is a straightforward application of Markov's Inequality and the Lemma of Borel--Cantelli.
\begin{Lemma}\label{lemma:upperbound}
$$
\limsup_{m\to\infty} |\mathcal{H}_m|^{1/m} \leq \frac{1}{z^\ast}\
\textrm{ almost surely.}
$$ 
\end{Lemma}
Eventually, we obtain the desired upper box-counting dimension.
\begin{Prop}\label{prop:upperbound}
$$
\overline{\mathrm{BD}}(\Lambda) \leq \frac{\log z^\ast}{\log \alpha}
$$
\end{Prop}
\begin{proof}
Denote by $N(\alpha^m)$ the number of balls of radius of at most $\alpha^m$ needed to
cover $\Lambda$. Then, for any $\varepsilon>0$, $N(\alpha^m)\leq
|\mathcal{H}_m|\leq \bigl(\frac{1}{z^\ast}+\varepsilon\bigr)^m$ almost surely
for sufficiently large $m$. Therefore,
$$
\overline{\mathrm{BD}}(\Lambda) = \limsup_{m\to\infty} -\frac{\log
  N(\alpha^m)}{\log \alpha^m} \leq
\limsup_{m\to\infty} -\frac{\log \bigl(\frac{1}{z^\ast}+\varepsilon\bigr)^m}{\log \alpha^m}=
-\frac{\log \bigl(\frac{1}{z^\ast}+\varepsilon\bigr)}{\log \alpha}.
$$
Letting  $\varepsilon\to 0$ proves the claim.
\end{proof}

\subsubsection{Lower Bound for Hausdorff Dimension}
\label{subsec:lower-bound}
In this section we will show that $\log z^\ast/\log\alpha$ is also a lower bound for the Hausdorff dimension of $\Lambda$.
From this we may then conclude existence of the box-counting dimension since
$\mathrm{HD}(\Lambda)\leq \underline{\mathrm{BD}}(\Lambda)\leq
\overline{\mathrm{BD}}(\Lambda).$ The main idea of the proof follows
\cite{lalley2000}. This idea\footnote{In this section the parameter $r$ is not identified with
$|\mathcal{I}|$  but is used as a parameter of the Galton--Watson trees
$\tau_{r}$ as in \cite{lalley2000}.} is to construct a sequence of embedded
Galton--Watson trees $\tau_{r}$ in the BRW such that the limit set
$\Lambda_{\tau_{r}}$ of the Galton--Watson trees are subsets of the limit set
$\Lambda$, see Section 6.3 in \cite{lalley2000}. As $r$ goes to infinity we will
have that $HD(\Lambda_{\tau_{r}})\to HD(\lambda)$.  
This approximation property relies mainly on the facts that particles travel essentially along geodesics segments and that limit sets of multi-type Galton--Watson trees are well understood. Both facts hold still true for free products of finite groups and the proof of the lower bound is analogous to  the one for free groups in \cite{lalley2000}, albeit technically more involved. The case of infinite factors need some extra care, since in this case particles do not necessarily travel along geodesics and infinite-type Galton--Watson processes are not so easy to handle. 
 To bypass these difficulties we approximate the infinite factors by increasing sequence of finite subgraphs. These subgraphs $\X_i^{(d)}$ are the subgraphs induced by the balls $B_i(d):=\{y\in\Gamma_i \mid l(y)\leq d\}, d\geq 1$.  Letting $d\to\infty$  will give the optimal bound $\log
z^\ast / \log \alpha$. 
\par
We add an
additional vertex $\dag$ to  $\X_i^{(d)}$, the ``tomb'', such that all edges in $\mathcal{X}_i$ exiting
$B_i(d)$ now lead to the tomb. The random walk $\bigl(Y^{(i,d)}_n\bigr)_{n\in\N_0}$ on $\X_i^{(d)}$ behaves like the random walk on
$\Gamma_i$, with the exception that a particle leaving 
$B_i(d)$ dies.  
We now build
the free product $\mathcal{X}^{(d)}$ from  the $\X_i^{(d)}$, whose
vertices are given by the set
$$
\Bigl\lbrace x_1\dots x_n\in\Gamma \,\Bigl|\, n\in\N, x_j\in\bigcup_{i\in\calI}
\X_i^{(d)}\setminus\{e_i,\dag\}, x_j\in \X_i^{(d)} \Rightarrow  x_{j+1}\notin
\X_i^{(d)} \Bigr\rbrace \cup \bigl\lbrace e,\dag\bigr\rbrace,
$$
where $\dag$ symbolizes the tomb. 
We identify  $x\in\X^{(d)}$ with the corresponding element in $\Gamma$. Analogously to Subsection
\ref{subsec:fp-rw}, we lift the random walks on the graphs $\X_i^{(d)}$ to a
random walk $\bigl( X_n^{(d)}\bigr)_{n\in\N_0}$ on $\X^{(d)}$ and define the associated BRW. 
We use the same notation (for Green functions, generating functions, etc.) as for the
random walk on $\Gamma$ itself but for reason of distinguishing we add
superscripts ``$(d)$'', that is, we write, e.g., $G^{(d)}(x,y|z)$ for the
corresponding Green function of the random walk on $\X^{(d)}$. All involved
generating functions on $\mathcal{X}^{(d)}$ have radii of convergence of at
least $R$.
\par
For any $x,y\in\Gamma$, we define $\overline{x:y}$ to be the set of
vertices $w\in\Gamma$ such that there is a geodesic from $x$ to $y$ which
passes through $w$. For $u\in\Gamma$, $d(u,\overline{x:y})$ is defined as the minimal
distance w.r.t. the  graph metric of $u$ to any element of $\overline{x:y}$. In the case of the coloured branching random walk on
$\X^{(d)}$, let $Z_{\infty}^{(d)}(y|x)$ be the overall number of blue particles arriving and
freezing at $y\in\mathcal{X}^{(d)}$ under the assumption that the branching random walk is
started with one blue particle at $x$. For $r\in\N$, we write $Z_{\infty,r}^{(d)}(y|x)$ for the overall number
of particles counted in $Z_{\infty}^{(d)}(y|x)$ whose trail remain within distance $r$
to a geodesic from $x$ to $y$. In other words, in all sites $u$ with
$d(u,\overline{x:y})>r$ every blue particle is coloured red. In the following we set $x_0:=x^{-1}_{1}$ for any $x=x_1\dots x_m\in\X^{(d)}$. The proofs of the two following lemmas are similar to the ones of Lemma 4 and Proposition 7 in \cite{lalley2000} and are therefore omitted.\footnote{The reader may find all the details in the \emph{arxiv.org} version of this paper.}
\begin{Lemma}\label{lemma:lim-inf-quotient} 
$$
\lim_{r\to\infty} \inf_{x=x_1\dots x_m\in \mathcal{X}^{(d)}} \biggl(\frac{
\prod_{j=1}^m \mathbb{E}Z_{\infty,r}^{(d)}(x_1\dots x_j|x_1\dots x_{j-1})}{\mathbb{E}Z_{\infty}^{(d)}(x|e)} \biggr)^{1/l(x)}=1.
$$
\end{Lemma}
For $x\in\mathcal{X}^{(d)}$, we define the event $E^{(d)}(x)$ that among the
particles counted in $Z_{\infty}^{(d)}(x|e)$ there is at least one particle whose trail
has not entered $\Gamma_1^\times$ and enters the set
$$
\bigl\lbrace y\in\mathcal{X}^{(d)} \mid l(y)=l(x) \bigr\rbrace
$$
first at $x$. Obviously, $Z_{\infty}^{(d)}(x|e)\geq 1$ on the event
$E^{(d)}(x)$ and hence $\Prob\bigl[E^{(d)}(x)\bigr]\leq \mathbb{E}
Z_{\infty}^{(d)}(x|e)$. 
\begin{Lemma}\label{lemma:e(x)}
$$
\lim_{k\to\infty} \biggl(\min_{{x=x_1\dots x_m \in\mathcal{X}^{(d)}:}
  \atop{m\in\N, x_1\notin\Gamma_1, l(x)=k}}
\frac{\Prob[E^{(d)}(x)]}{\mathbb{E} Z^{(d)}_{\infty}(x|e)}\biggr)^{1/k} =1.
$$
\end{Lemma}
Analogously to (\ref{equ:scriptFi+}) and (\ref{equ:scriptFi}), we define for
$i\in\mathcal{I}$ and $d\in\N$
\begin{eqnarray}
\mathcal{L}_i^{(d)+}(\lambda|z) & := & \sum_{x\in\Gamma_i^\times} L^{(d)}(e,x|\lambda)\,z^{l(x)}=
\sum_{x\in\Gamma_i^\times} L^{(d)}_i\bigl(e_i,x\bigl|\xi_i^{(d)}(\lambda)\bigr)\,z^{l(x)},\notag \\
\mathcal{L}^{(d)}_i(\lambda|z) & := & \sum_{n\geq 1}\sum_{\substack{x=x_1\dots x_n
  \in\X^{(d)}:\\ \tau(x_1)=i}} L^{(d)}(e,x|\lambda)\,z^{l(x)}\notag \\
  & = & 
\mathcal{L}_i^{(d)+}(\lambda|z)\Bigl(1+\sum_{j\in\mathcal{I}\setminus\{i\}}
\mathcal{L}^{(d)}_j(\lambda|z)\Bigr). \label{eq:Lid}
\end{eqnarray}

Writing 
$\mathcal{L}^{(d)}(\lambda|z):=1+\sum_{i\in\mathcal{I}}\mathcal{L}^{(d)}_i(\lambda|z)$  we get analogously to Equation (\ref{equ:F-equation}): 
$$
\mathcal{L}^{(d)}(\lambda|z)=\frac{1}{1-\sum_{i\in\mathcal{I}} \frac{\mathcal{L}_i^{(d)+}(\lambda|z)}{1+\mathcal{L}_i^{(d)+}(\lambda|z)}}.
$$
Since every function $\mathcal{L}_i^{(d)+}(\lambda|z)$ is convergent and
strictly increasing for all
$z\geq 0$ there is some unique $z_{d,\mathcal{L}}^\ast>0$ such that
$\sum_{i\in\mathcal{I}}
\mathcal{L}_i^{(d)+}(\lambda|z_{d,\mathcal{L}}^\ast)/\bigl(1+\mathcal{L}_i^{(d)+}(\lambda|z_{d,\mathcal{L}}^\ast)\bigr)=1$. The radius of
convergence of $\mathcal{L}^{(d)}(\lambda|z)$ is then given by $z_{d,\mathcal{L}}^\ast$.
\par
We define for $k\in\N$
$$
S^\ast_k := \bigl\lbrace x_1\dots x_s\in \X^{(d)} \,\bigl|\, s\in\N,l(x)=k, x_1\notin\Gamma_1,x_s\in\Gamma_1 \bigr\rbrace.
$$
Since we excluded the case $|\mathcal{I}|=2=|\Gamma_1|=|\Gamma_2|$ we have  that $S^\ast_2\neq \emptyset$ and $S^\ast_3\neq \emptyset$. Therefore, $S^\ast_k\neq
\emptyset$ for all $2\leq k\in\N$.
\begin{Lemma}\label{lemma:sum-root}
$$
\limsup_{k\to\infty}\biggl( \sum_{x\in S^\ast_k}
\Prob[E^{(d)}(x)]\biggr)^{1/k} =  \frac{1}{z_{d,\mathcal{L}}^\ast}.
$$
\end{Lemma}
\begin{proof}
By Lemma \ref{lemma:e(x)}, we have $\Prob[E^{(d)}(x)]\geq (1-\varepsilon)^k
\mathbb{E} Z_\infty^{(d)}(x|e)$ uniformly for all $x$ with $l(x)=k$ if $k$ is large enough. Recall also $\Prob[E^{(d)}(x)]\leq \mathbb{E} Z_\infty^{(d)}(x|e)$.
Thus, it is sufficient to prove
$$
\limsup_{k\to\infty}\biggl( \sum_{x\in S^\ast_k}
\mathbb{E} Z_\infty^{(d)}(x|e) \biggr)^{1/k} =
\frac{1}{z_{d,\mathcal{L}}^\ast}. 
$$
Since
$$
\sum_{x\in S^\ast_k}
\mathbb{E} Z_\infty^{(d)}(x|e)  = 
\sum_{x\in S^\ast_k}
F^{(d)}(e,x|\lambda) 
=
\sum_{x\in S^\ast_k}
\frac{G^{(d)}(e,e|\lambda)}{G^{(d)}(x,x|\lambda)} L^{(d)}(e,x|\lambda)
$$
and $1\leq G^{(d)}(x,x|\lambda) \leq G(x,x|\lambda) = G(e,e|\lambda)<\infty$ we have
\begin{equation}\label{equ:Ld-conv}
\limsup_{k\to\infty}\biggl(\sum_{x\in S^\ast_k} L^{(d)}(e,x|\lambda)\biggr)^{1/k} = \limsup_{k\to\infty}\Bigl( \sum_{x\in S^\ast_k} \mathbb{E}
Z_\infty^{(d)}(x|e) \Bigr)^{1/k}.
\end{equation}
To determine the left-hand side of   (\ref{equ:Ld-conv}) 
we define further generating functions:
\begin{eqnarray*}
\mathcal{L}^{(d)}_{\neg 1,1}(\lambda|z) & := & \sum_{n\geq 2} \sum_{\substack{x=x_1\dots x_n\in\mathcal{X}^{(d)}:\\ x_1\notin\Gamma_1^\times,x_n\in\Gamma_1^\times}} L^{(d)}(e,x|\lambda)\,z^{l(x)},\\
\mathcal{L}^{(d)*}(\lambda|z) & := & \sum_{n\geq 1} \sum_{\substack{x=x_1\dots x_n\in\mathcal{X}^{(d)}:\\ x_1,x_2,\dots,x_n\notin\Gamma_1}} L^{(d)}(e,x|\lambda)\,z^{l(x)}.
\end{eqnarray*}
For $k\in\N$, the coefficient of $z^k$ in $\mathcal{L}^{(d)}_{\neg 1,1}(\lambda|z)$ is just $\sum_{x\in S^\ast_k} L^{(d)}(e,x|\lambda)$. Due to Equation (\ref{eq:Lid}) we have
$$
\mathcal{L}_1^{(d)}(\lambda|z)=\mathcal{L}_1^{(d)+}(\lambda|z)\cdot \Bigr(1+
\sum_{i\in\mathcal{I}\setminus\{1\}}\mathcal{L}_i^{(d)}(\lambda|z)\Bigr),
$$ and hence
the function $\mathcal{L}^{(d)}_{\neg 1}(\lambda|z):=
1+\sum_{i\in\mathcal{I}\setminus\{1\}}\mathcal{L}_i^{(d)}(\lambda|z)$
must have the same radius of convergence as $\mathcal{L}^{(d)}(\lambda|z)$,
which is $z_{d,\mathcal{L}}^\ast$. Moreover, we have the following relations:
\begin{eqnarray*}
\mathcal{L}^{(d)}_{\neg 1}(\lambda|z) & = & 1+\mathcal{L}^{(d)}_{\neg 1,1}(\lambda|z)\bigl(1+\mathcal{L}^{(d)*}(\lambda|z)\bigr) + \mathcal{L}^{(d)*}(\lambda|z),\\
\mathcal{L}^{(d)}_{\neg 1,1}(\lambda|z)& \geq & \mathcal{L}^{(d)*}(\lambda|z)\cdot  \mathcal{L}_1^{(d)+}(\lambda|z).
\end{eqnarray*}
Since $\mathcal{L}^{(d)}_{\neg 1,1}(\lambda|z),\mathcal{L}^{(d)*}(\lambda|z)\leq \mathcal{L}^{(d)}_{\neg 1}(\lambda|z)$, the function $\mathcal{L}^{(d)}_{\neg 1,1}$ has also radius of convergence of $z_{d,\mathcal{L}}^\ast$. 
\end{proof}
Now we show that $z_{d,\calL}^\ast$ tends to $z^\ast$ as $d\to\infty$. Since
$z_{d,\calL}^\ast$ is strictly decreasing as $d$ grows and due to
\begin{equation}\label{equ:Ld-conv2}
\lim_{d\to\infty} L^{(d)}(e,x|\lambda) = L(e,x|\lambda)=F(e,x|\lambda)
\end{equation}
we have $z_\infty=\lim_{d\to\infty} z_{d,\calL}^\ast\geq z^\ast$.
Assume now for a moment that  $z^\ast<z_\infty$. Then $\mathcal{F}^+_i(\lambda|z_\infty)<\infty$:
 indeed, assume that $\lim_{d\to\infty}\mathcal{L}_j^{(d)+}(\lambda|z_\infty)=\mathcal{F}^+_j(\lambda|z_\infty)=\infty$
for some $j\in\mathcal{I}$. Then we get the following contradiction:
\begin{equation}\label{equ:F<infty}
1=\lim_{d\to\infty} \sum_{i\in\mathcal{I}}
\frac{\mathcal{L}_i^{(d)+}(\lambda|z_{d,\calL}^\ast)}{1+\mathcal{L}_i^{(d)+}(\lambda|z_{d,\calL}^\ast)}
\geq
\lim_{d\to\infty} \sum_{i\in\mathcal{I}}
\frac{\mathcal{L}_i^{(d)+}(\lambda|z_\infty)}{1+\mathcal{L}_i^{(d)+}(\lambda|z_\infty)}>1,
\end{equation}
since
$\mathcal{L}_j^{(d)+}(\lambda|z_\infty)/\bigl(1+\mathcal{L}_j^{(d)+}(\lambda|z_\infty)\bigr)$
is arbitrarily close to $1$ if $d$ is large enough. Hence,
$\mathcal{F}^+_i(\lambda|z_\infty)<\infty$. Now $z_\infty >z^\ast$
yields the following contradiction:
$$
1= \lim_{d\to\infty} \sum_{i\in\mathcal{I}}
 \frac{\mathcal{L}^{(d)+}_i(\lambda|z_{d,\calL}^\ast)}{1+\mathcal{L}^{(d)+}_i(\lambda|z_{d,\calL}^\ast)}
\geq \limsup_{d\to\infty} \sum_{i\in\mathcal{I}}
 \frac{\mathcal{L}^{(d)+}_i(\lambda|z_{\infty})}{1+\mathcal{L}^{(d)+}_i(\lambda|z_{\infty})}
= \sum_{i\in\mathcal{I}}
 \frac{\mathcal{F}^{+}_i(\lambda|z_{\infty})}{1+\mathcal{F}^{+}_i(\lambda|z_{\infty})}>1.
$$
which produces  a contradiction. Thus,
\begin{equation}\label{equ:approx}
\lim_{d\to\infty} z_{d,\mathcal{L}}^\ast=z^\ast.
\end{equation}
Let $2\leq k\in\N$ arbitrary, but fixed.
Similar to \cite{lalley2000} we define an  embedded Galton--Watson process of the BRW on the free product
$\X^{(d)}$. For $n\in\N_0$, we define generations
$\mathrm{gen}(n)$  $S^\ast_{nk}$ and distinguished  particles $\zeta_{x}$ associated to vertices $x\in \mathrm{gen}(n)$ inductively as follows:
\begin{enumerate}
\item $\mathrm{gen}(0):=\{e\}$ consists of one particle $\zeta_e$ located at $e$.
\item 
   $y\in S^\ast_{(n+1)k}$ belongs to
  $\mathrm{gen}(n+1)$ if and only if there exists a distinguished particle $\zeta_{x}$ in $\mathrm{gen}(n)$ such that some of its offspring particles counted in $Z_\infty^{(d)}(y|x)$ has a trail which
\begin{enumerate}
\item remains in the set
$$
\Gamma(x):=\{y\in\Gamma \mid y \textrm{ has the form } xw_1\dots w_s \textrm{ with }
  w_1\notin \Gamma_1, s\geq 1\}\cup\{x\},
$$
\item hits the set $\bigl\lbrace w\in\X^{(d)} \,\bigl|\, l(w)=(n+1)k\bigr\rbrace$ first at $y$.
\end{enumerate}
\item The first particle hitting $y\in S^\ast_{(n+1)k}$ becomes the distinguished particle $\zeta_y$.
\end{enumerate}
%
Let $\phi_{n}$ denote the number of particles in generation $n$. Since we have the same offspring distribution at every $x\in
S^\ast_{nk}$, $(\phi_{n})_{n\geq 0}$ defines a Galton--Watson process with mean
$M_{d,k}$.
\begin{Cor}
$$
\limsup_{k\to\infty} M_{d,k}^{1/k} = \frac{1}{z_{d,\mathcal{L}}^\ast}.
$$
\end{Cor}
\begin{proof}
The claim follows directly with Lemma \ref{lemma:sum-root} since $M_{d,k} = \sum_{x\in S^\ast_k} \Prob[E^{(d)}(x)]$.
%
\end{proof}
Applying Hawkes' Theorem as in Corollary 7 in \cite{lalley2000} together with Equation
(\ref{equ:approx}) yields the following Corollary.

\begin{Cor}\label{cor:boxdim-lambda}
$ \mathrm{HD}(\Lambda \cap \Omega_\infty) \geq \frac{\log z^\ast}{\log\alpha}$.
\end{Cor}
\begin{proof}[Proof of Theorem \ref{thm:boxdim}.]
The following chains of inequalities summarize the previous results and finish the proof of the theorem:
$$
\begin{array}{c}
\frac{\log z^\ast}{\log\alpha} \leq \mathrm{HD}(\Lambda) \leq
\underline{\mathrm{BD}}(\Lambda) \leq
\overline{\mathrm{BD}}(\Lambda)\leq \frac{\log z^\ast}{\log\alpha},\\[2ex]
\frac{\log z^\ast}{\log\alpha} \leq \mathrm{HD}(\Lambda \cap \Omega_\infty) \leq
\underline{\mathrm{BD}}(\Lambda \cap \Omega_\infty) \leq
\overline{\mathrm{BD}}(\Lambda \cap \Omega_\infty)\leq 
\overline{\mathrm{BD}}(\Lambda)\leq
\frac{\log z^\ast}{\log\alpha}.
\end{array}
$$
\end{proof}
\begin{proof}[Proof of Corollary \ref{cor:omega-i-do-not-contribute}.]
It is well-known that the Hausdorff dimension of a countable union
$\bigcup_{i} B_i$ of sets $B_i\subseteq \Omega$
equals the supremum of the Hausdorff dimensions of the single sets $B_i$. Thus,
$$
\mathrm{HD}(\Lambda \cap \Omega_i) =
\sup_{x\in \Gamma: \tau(x)\neq i}  \mathrm{HD}(\Lambda \cap
x\Omega_i^{(0)}) \leq
\sup_{x\in \Gamma: \tau(x)\neq i}  \overline{\mathrm{BD}}(\Lambda \cap
x\Omega_i^{(0)}).
$$
For arbitrary, but fixed $x\in \Gamma$ with $\tau(x)\neq i$, denote by
$\mathcal{H}^{(x)}_m$ the vertices $y\in x\Gamma_i$ with $l(y)=l(x)+m$, which are visited by the
branching random walk. Therefore,
$$
\mathbb{E}|\mathcal{H}^{(x)}_m| \leq \sum_{y\in\Gamma_i: l(y)=m} F(e,xy|\lambda) =
F(e,x|\lambda) \sum_{y\in\Gamma_i: l(y)=m} F(e,y|\lambda).
$$
Define
$$
\mathcal{F}_{x,i}^+(\lambda|z):= F(e,x|\lambda) \sum_{m\geq 1} \sum_{y\in\Gamma_i: l(y)=m} F(e,y|\lambda)\,z^m.
$$
The radius of convergence of $\mathcal{F}_{x,i}^+(\lambda|z)$ is obviously $R(\mathcal{F}_i^+)$. Therefore, Lemma \ref{lemma:conv-radius2} yields $\limsup_{m\to\infty}\bigl(\mathbb{E}|\mathcal{H}^{(x)}_m|\bigr)^{1/m}\leq
1/R(\mathcal{F}_{i}^+)<1/z^\ast$. The rest follows analogously to the proofs of Lemma \ref{lemma:upperbound} and Proposition \ref{prop:upperbound}.
\end{proof}

\subsection{Proof of Theorem \ref{thm:boxdim-omega} and Corollary \ref{cor:boxdim-omega}}

In order to prove Theorem \ref{thm:boxdim-omega} we can follow the argumentation of the proof of Theorem \ref{thm:boxdim}. For this purpose, we
define for $m\in\N$ and $i\in\mathcal{I}$
$$
\begin{array}{c}
S_i(m):=|\{x\in\Gamma_i \mid l(x)=m\}|, \quad  S(m) :=  |\{x\in\Gamma \mid l(x)=m\}|,\\[2ex]
S^{(i)}(m)  :=  |\{x=x_1\dots x_n\in S(m) \mid n\in\N, x_1\in\Gamma_i\}|.
\end{array}
$$
To cover $\Omega$ by balls of radius $\alpha^m$ we need at least
$S(m-1)$ balls: indeed, for all $x,y\in\Gamma$, $x\neq y$, with $l(x)=l(y)=m-1$ we can choose $v_x\in\Gamma^\times_\ast\setminus \Gamma_{\tau(x)}$ and $v_y\in\Gamma^\times_\ast\setminus \Gamma_{\tau(y)}$; then all balls of the form $B(\omega_1,\alpha^m)$ and $B(\omega_2,\alpha^m)$, where
$xv_x$ lies in the
$\omega_1$-component of $\mathcal{X}\setminus B_{m-1}$ and $yv_y$ in the
$\omega_2$-component, do not intersect. Apparently, we need at most $S(m)$ balls
of radius $\alpha^m$ to cover $\Omega$. Obviously, the same holds for covering $\Omega_\infty$. We are now interested in the
behaviour of $S(m)^{1/m}$ as $m\to\infty$. We define
\begin{eqnarray*}
\mathcal{S}_i^+(z) & := & \sum_{m\geq 1} S_i(m)\,z^m,\quad
\mathcal{S}_i(z) := \sum_{m\geq 1} S^{(i)}(m)\,z^m,\\
\mathcal{S}(z) & := & \sum_{m\geq 0} S(m)\,z^m=1+\sum_{i\in\mathcal{I}} \mathcal{S}^{(i)}(z).
\end{eqnarray*}
Analogously to the computations in Section \ref{subsec:upperbound} -- we just
replace the functions $\mathcal{F}_i^+(\lambda|z)$, $\mathcal{F}_i(\lambda|z)$
and $\mathcal{F}(\lambda|z)$ by the functions $\mathcal{S}_i^+(z)$,
$\mathcal{S}_i(z)$ and $\mathcal{S}(z)$ -- we get
\begin{equation}\label{equ:S-formula}
\mathcal{S}(z) = \frac{1}{1-\sum_{i\in\mathcal{I}} \frac{\mathcal{S}_i^+(z)}{1+\mathcal{S}_i^+(z)}}.
\end{equation}
\begin{Lemma}
$$
\lim_{m\to\infty} S(m)^{1/m}=\frac{1}{z^\ast_{\mathcal{S}}}<1,
$$
where $z_{\mathcal{S}}^\ast$ is the smallest positive real number with 
$$
\sum_{i\in\mathcal{I}}
\frac{\mathcal{S}_i^+(z_{\mathcal{S}}^\ast)}{1+\mathcal{S}_i^+(z_{\mathcal{S}}^\ast)}=1.
$$
\end{Lemma}
\begin{proof}
Obviously, $R(\mathcal{S})\leq R(\mathcal{F})<1$ since $F(e,x|\lambda)<1$ for all $x\in\Gamma\setminus\{e\}$. The equation $R(\mathcal{S})=z_{\mathcal{S}}^\ast$ follows now analogously to the proof of Lemma \ref{lemma:conv-radius2}. This yields
$$
\limsup_{m\to\infty}S(m)^{1/m}=\frac{1}{z^\ast_{\mathcal{S}}}=\frac{1}{R(\mathcal{S})}> 1.
$$ 
Thus, it is sufficient to prove convergence of $S(m)^{1/m}$ as $m\to\infty$. By
transitivity of $\Gamma$, we have $S(m)S(n)\geq S(m+n)$ for all
$m,n\in\N$. Therefore,  $\log S(m) + \log S(n)\geq \log S(m+n)$, that is,
$\bigl(\log S(m)\bigr)_{m\in\N}$ forms a subadditive sequence. By Fekete's Lemma, $\frac{1}{m} \log S(m)=\log S(m)^{1/m}$ converges
to some constant $s$, that is, $S(m)^{1/m}$ converges to $e^s$, which must
equal $1/z^\ast_{\mathcal{S}}$. 
\end{proof}
\begin{Rem}\label{remark:z-ast}
One can show analogously to Lemma \ref{lemma:conv-radius2} that
$z_{\mathcal{S}}^\ast<R(\mathcal{S}_i^+)$, where $R(\mathcal{S}_i^+)$ is the
radius of convergence of $\mathcal{S}_i^+(z)$. In particular,
$z_{\mathcal{S}}^\ast$ is the radius of convergence of $\mathcal{S}(z)$.
\end{Rem}
We can conclude by giving a formula for $\mathrm{BD}(\Omega)$ and
observing that the box-counting dimension of $\Omega$ results from the
dimension of $\Omega_\infty$.
\begin{Prop}\label{prop:boxdim-omega}
$$
\mathrm{BD}(\Omega)=\mathrm{BD}(\Omega_\infty) = \frac{\log
  z^\ast_{\mathcal{S}}}{\log\alpha}.
$$
\end{Prop}
\begin{proof}
Recall the remarks at the beginning of this section concerning the minimal
and maximal number of balls needed to cover $\Omega_\infty$. This yields
\begin{eqnarray*}
\underline{\mathrm{BD}}(\Omega) &\geq &
\underline{\mathrm{BD}}(\Omega_\infty )\geq
\liminf_{m\to\infty} -\frac{\log
  S(m-1)}{\log \alpha^m} \\
 &= & \liminf_{m\to\infty} -\frac{\log
  S(m-1)^{1/(m-1)}}{\log \alpha}\frac{m-1}{m}=\frac{\log z^\ast_{\mathcal{S}}}{\log\alpha}.
\end{eqnarray*}
Analogously,
$$
\overline{\mathrm{BD}}(\Omega_\infty )\leq \overline{\mathrm{BD}}(\Omega)\leq \limsup_{m\to\infty} -\frac{\log
  S(m)}{\log \alpha^m} = \limsup_{m\to\infty} -\frac{\log
  S(m)^{1/m}}{\log \alpha}=\frac{\log z^\ast_{\mathcal{S}}}{\log\alpha}.
$$
Both inequality chains together yield the formula for the box-counting
dimension. 
\end{proof}
Finally, we can prove the formula for the Hausdorff dimensions of $\Omega$ and $\Omega_\infty$.
\begin{proof}[Proof of Theorem \ref{thm:boxdim-omega}.]
It is sufficient to show that $\mathrm{HD}(\Omega_\infty)\geq \frac{\log
  z^\ast_\mathcal{S}}{\log\alpha}$. Define for $d,k\in\N$ and $i\in\mathcal{I}$
\begin{eqnarray*}
S_i^{(d)+}(k) & = & \bigl|\bigr\lbrace x\in\mathcal{X}_i^{(d)}\setminus \{\dag\} \mid
l(x)=k\bigr\rbrace\bigl|,\
S^{(d)}(k) =  \bigl|\bigr\lbrace x\in\mathcal{X}^{(d)}\setminus \{\dag\} \mid
l(x)=k\bigr\rbrace\bigl|,\\
S^{(d)}_i(k)& = & \bigl|\bigr\lbrace x_1\dots x_s\in\mathcal{X}^{(d)}\setminus \{\dag\} \mid s\in\N,
l(x)=k, x_1\in\Gamma_i\bigr\rbrace\bigl|,\\
S^{(d)}_{\neg 1}(k)& = & \bigl|\bigr\lbrace x_1\dots x_s\in\mathcal{X}^{(d)}\setminus \{\dag\}
\mid s\in\N, l(x)=k, x_1\notin\Gamma_1 \bigr\rbrace\bigl|,\\
S^{(d)}_{\neg 1,\neg 1}(k)& = & \bigl|\bigr\lbrace x_1\dots x_s\in\mathcal{X}^{(d)}\setminus \{\dag\}
\mid s\in\N, l(x)=k, x_1,x_s\notin\Gamma_1 \bigr\rbrace\bigl|,\\
S^{(d)}_{\neg 1,1}(k)& = & \bigl|\bigr\lbrace x_1\dots x_s\in\mathcal{X}^{(d)}\setminus \{\dag\}
\mid s\in\N, l(x)=k, x_1\notin\Gamma_1, x_s\in\Gamma_1 \bigr\rbrace\bigl|.
\end{eqnarray*}
The associated generating functions are given by
\begin{eqnarray*}
\mathcal{S}_i^{(d)+}(z) &=& \sum_{k\geq 1} S_i^{(d)+}(k)\,z^k,
\quad \mathcal{S}^{(d)}(z)=\sum_{k\geq 0} S^{(d)}(k)\,z^k,\\
\mathcal{S}^{(d)}_i(z)&=&\sum_{k\geq 1} S^{(d)}_i(k)\,z^k,\quad 
\mathcal{S}^{(d)}_{\neg 1}(z) = \sum_{k\geq 1} S^{(d)}_{\neg 1}(k)\,z^k,\\
\mathcal{S}^{(d)}_{\neg 1,\neg 1}(z)&=&\sum_{k\geq 1} S^{(d)}_{\neg 1,\neg 1}(k)\,z^k,\quad
\mathcal{S}^{(d)}_{\neg 1,1}(z) = \sum_{k\geq 1} S^{(d)}_{\neg 1,1}(k)\,z^k.
\end{eqnarray*}
Once again we can write
$$
\mathcal{S}^{(d)}(z)=\frac{1}{1-\sum_{i\in\mathcal{I}} \frac{\mathcal{S}^{(d)+}_i(z)}{1+\mathcal{S}^{(d)+}_i(z)}}
$$
and obtain
$$
\mathcal{S}^{(d)+}_2(z) \mathcal{S}^{(d)+}_1(z) \mathcal{S}^{(d)}_{\neg 1,\neg
  1}(z) \mathcal{S}^{(d)+}_1(z) \leq \mathcal{S}^{(d)}_{\neg 1,1}(z) = 
\mathcal{S}^{(d)}_{\neg 1}(z) - \mathcal{S}^{(d)}_{\neg 1,\neg 1}(z).
$$
Thus, $\mathcal{S}^{(d)}_{\neg 1,1}(z)$ and $\mathcal{S}^{(d)}_{\neg 1}(z)$ have the
same radius of convergence. Moreover,
$$
\mathcal{S}^{(d)+}_2(z) \mathcal{S}^{(d)}_1(z) \leq \mathcal{S}^{(d)}_{\neg
  1}(z)
= \mathcal{S}^{(d)}(z)-\mathcal{S}^{(d)}_1(z)-1.
$$
That is, $\mathcal{S}^{(d)}_{\neg 1,1}(z)$ and $\mathcal{S}^{(d)}(z)$ have the same
radius of convergence, which is given by $z^\ast_{d,\mathcal{S}}$, the smallest positive solution satisfying
$$
1=\sum_{i\in\mathcal{I}} \frac{\mathcal{S}^{(d)+}_i(z)}{1+\mathcal{S}^{(d)+}_i(z)}.
$$
Since  $z^\ast_{d,\mathcal{S}}$ is strictly decreasing as $d\to\infty$ we have that $\lim_{d\to\infty} z^\ast_{d,\mathcal{S}}=z^\ast_\mathcal{S}$. This can be seen by contradiction.  Indeed, if
$\lim_{d\to\infty}
z^\ast_{d,\mathcal{S}}=z^\ast_{\infty,\mathcal{S}}>z^\ast_{\mathcal{S}}$  then
$\mathcal{S}^+_i(z^\ast_{\infty,\mathcal{S}})<\infty$ for all $i\in\calI$
(this is proven analogously as explained in Equation (\ref{equ:F<infty})) 
and therefore
$$
1=\lim_{d\to\infty} \sum_{i\in\mathcal{I}}
\frac{\mathcal{S}^{(d)+}_i(z^\ast_{d,\mathcal{S}})}{1+\mathcal{S}^{(d)+}_i(z^\ast_{d,\mathcal{S}})}
\geq
\lim_{d\to\infty} \sum_{i\in\mathcal{I}}
\frac{\mathcal{S}^{(d)+}_i(z^\ast_{\infty,\mathcal{S}})}{1+\mathcal{S}^{(d)+}_i(z^\ast_{\infty,\mathcal{S}})}=\sum_{i\in\mathcal{I}} \frac{\mathcal{S}^+_i(z^\ast_{\infty,\mathcal{S}})}{1+\mathcal{S}^+_i(z^\ast_{\infty,\mathcal{S}})}>1,
$$
a contradiction. Thus,
$$
\Bigl(S^{(d)}_{\neg 1,1}(k)\Bigr)^{1/k} \xrightarrow{k\to\infty} \frac{1}{z^\ast_{d,\mathcal{S}}}
\xrightarrow{d\to\infty} \frac{1}{z^\ast_{\mathcal{S}}}.
$$
We can embed a ``deterministic'' Galton--Watson tree into the free product
analogously to Subsection \ref{subsec:lower-bound}, where each generation has
exactly $S^{(d)}_{\neg 1,1}(k)$ descendants. By Hawkes's Theorem, the Hausdorff dimension
of the boundary of the  embedded tree is bounded from below by $\log z^\ast_{d,\mathcal{S}}/\log \alpha$,
and therefore $\mathrm{HD}(\Omega_\infty)\geq \log z^\ast_{\mathcal{S}}/\log
\alpha$.
\end{proof}
\begin{proof}[Proof of Corollary \ref{cor:boxdim-omega}.]
Analogously to the proof of Corollary \ref{cor:omega-i-do-not-contribute} and
by Remark \ref{remark:z-ast}, we can use the property $\mathrm{HD}(\cup_{i} B_i)=\sup_i \mathrm{HD}(B_i)$ for all countable unions of sets $B_i\subseteq \Omega$ in order to show that
$$
\mathrm{HD}(\Omega_i)=
\sup_{x\in\Gamma: \tau(x)\neq i} \mathrm{HD}(x\Omega_i^{(0)}) \leq
\overline{\mathrm{BD}}(\Omega_i^{(0)})<\mathrm{BD}(\Omega_\infty)=\mathrm{HD}(\Omega_\infty).
$$
\end{proof}

\subsection{Proof of Theorem \ref{Th:lambda-function}.}
\label {sec:dim-comparision}
\begin{proof}[Proof of Theorem \ref{Th:lambda-function} (1).]
In the following we write $z^\ast=z^\ast(\lambda)$ in order to distinguish the
solutions of (\ref{equ:z-ast}) for different values of $\lambda$. Note that
$z^\ast(\lambda_1)>z^\ast(\lambda_2)$ if $\lambda_1<\lambda_2$. This implies
the strictly increasing behaviour of $\Phi$ in the interval $(1,R]$. Recall
that the BRW does almost surely \textit{not} survive in the limit case $\lambda = 1$, yielding $\Phi(1)=0$. Moreover, if $\lambda>R$ then the BRW is recurrent and thus $\mathrm{HD}(\Lambda)=\mathrm{HD}(\Omega)$.
\end{proof}
The proof of Theorem \ref{Th:lambda-function} (2) splits up into the following two lemmas:
\begin{Lemma}\label{lemma:phi-continuity}
$\Phi$ is continuous in $[1,\infty)\setminus\{R\}$ and continuous from the left  at $\lambda=R$.
\end{Lemma}
\begin{proof}
In order to prove continuity of $\Phi$, it is sufficient to
prove continuity of the mapping $\lambda \mapsto z^\ast=z^\ast(\lambda)$. 
First, we prove continuity from the left at $\lambda_0\in(1,\infty)$. For this purpose, let $(\lambda_n)_{n\in\N}$ be a sequence of strictly increasing real numbers
with $\lambda_n<\lambda_0$ and $\lim_{n\to\infty} \lambda_n=\lambda_0$. We use a proof by contradiction. Assume
$z_0:=\lim_{n\to\infty} z^\ast(\lambda_n) > z^\ast(\lambda_0)$ (by simple
domination arguments,
$z^\ast(\lambda_n)$ can not be less than $z^\ast(\lambda_0)$). We have that $z^\ast(\lambda_n)$ is strictly decreasing and 
$$ 
\mathcal{F}_i^+\bigr(\lambda_n\bigl|z^\ast(\lambda_0)\bigr) + \xi_i(1) \bigl(z_0-z^\ast(\lambda_0)\bigr) \leq 
\mathcal{F}_i^+\bigr(\lambda_n\bigl|z_0\bigr) <\infty. 
$$
Here we used the fact that the coefficient of $z$ in $\mathcal{F}_i^+(\lambda|z)$ is at least $\xi_i(1)$. We set $c:=\xi_i(1)
\bigl(z_0-z^\ast(\lambda_0)\bigr)$. Since $f(x)/\bigl(1+f(x)\bigr)$ is strictly increasing in $[1,\infty)$ if $f(x)$ is a strictly increasing function on $[1,\infty)$  we get the following contradiction:
\begin{eqnarray*}
1  &= &  \lim_{n\to\infty} \sum_{i\in\calI}
\frac{\mathcal{F}_i^+\bigr(\lambda_n\bigl|z^\ast(\lambda_n)\bigr)}{1+\mathcal{F}_i^+\bigr(\lambda_n\bigl|z^\ast(\lambda_n)\bigr)}\\
&\geq &  \limsup_{n\to\infty} \sum_{i\in\calI}
\frac{\mathcal{F}_i^+\bigr(\lambda_n\bigl|z_0\bigr)}{1+\mathcal{F}_i^+\bigr(\lambda_n\bigl|z_0\bigr)}\\
&\geq & \limsup_{n\to\infty} \sum_{i\in\calI}
\frac{\mathcal{F}_i^+\bigr(\lambda_n\bigl|z^\ast(\lambda_0)\bigr)+c}{1+\mathcal{F}_i^+\bigr(\lambda_n\bigl|z^\ast(\lambda_0)\bigr)+c}\\
& = & \sum_{i\in\calI}
\frac{\mathcal{F}_i^+\bigr(\lambda_0\bigl|z^\ast(\lambda_0)\bigr)+c}{1+\mathcal{F}_i^+\bigr(\lambda_0\bigl|z^\ast(\lambda_0)\bigr)+c}\\
&>& \sum_{i\in\calI}
\frac{\mathcal{F}_i^+\bigr(\lambda_0\bigl|z^\ast(\lambda_0)\bigr)}{1+\mathcal{F}_i^+\bigr(\lambda_0\bigl|z^\ast(\lambda_0)\bigr)} =1.
\end{eqnarray*} 
Thus, $\lim_{n\to\infty} z^\ast(\lambda_n) = z^\ast(\lambda_0)$.
\par
Since $\mathrm{HD}(\Lambda)=\mathrm{HD}(\Omega)$ for all $\lambda>R$, it remains to prove continuity from the right for $\lambda_0\in (1,R)$. We make a case distinction whether $\xi_i(\lambda_0) < 1$ or not. If $\xi_i(\lambda_0) < 1$ then $\mathcal{F}_i^+(\lambda_0+\delta|1)<\infty$ for all $\delta>0$ with $\xi_i(\lambda_0+\delta)<1$ according to (\ref{equ:fi+}). Moreover, $z^\ast(\lambda_0)<1$. Therefore, continuity from the right follows directly from the Implicit Function Theorem, since $z^\ast=z^\ast(\lambda)$ is given by the equation 
$$
1=\sum_{i\in\mathcal{I}} \frac{\mathcal{F}_i^+\bigl(\lambda\mid z^\ast(\lambda)\bigr)}{1+\mathcal{F}_i^+\bigl(\lambda\mid z^\ast(\lambda)\bigr)}.
$$
We note that the derivative $\partial\mathcal{F}_i^+(\lambda | z)/\partial z$ evaluated at $z=z^\ast(\lambda)$ is positive and finite, since $z^\ast(\lambda)$ is strictly smaller than the radius of convergence of $\mathcal{F}_i^+(\lambda | z)$; see Lemma \ref{lemma:conv-radius2}.
\par
Now we turn to the case $\xi_i(\lambda_0)\geq 1$.
Let $(\lambda_n)_{n\in\N}$ be a sequence of strictly decreasing real numbers
with $\lambda_0 < \lambda_n< R$ and $\lim_{n\to\infty} \lambda_n=\lambda_0$. Assume
$z_0:=\lim_{n\to\infty} z^\ast(\lambda_n) < z^\ast(\lambda_0)$ (by simple
domination arguments,
$z^\ast(\lambda_n)$ can not be larger than $z^\ast(\lambda_0)$). Observe that
$z^\ast(\lambda_n)$ is strictly increasing. By (\ref{equ:rF-rFi}), there is
$C:=\sqrt{1+\xi_2(1)/(2|\mathrm{supp}(\mu_1)|)}>1$ such that $C  z^\ast(\lambda_n)\leq R(\mathcal{F}_i^+)$ for all $n\in\N$. Choose $\tilde C\in(1,C)$ such that $\tilde C  z_0<z^\ast(\lambda_0)$ and choose $N\in\N$ large enough such that $\tilde C  z^\ast(\lambda_n)\geq z_0$ for all $n\geq N$. Therefore,
\begin{eqnarray*}
1 & = & \lim_{n\to\infty} \sum_{i\in\calI}
\frac{\mathcal{F}_i^+\bigr(\lambda_n\bigl|z^\ast(\lambda_n)\bigr)}{1+\mathcal{F}_i^+\bigr(\lambda_n\bigl|z^\ast(\lambda_n)\bigr)}\\
&\leq & \lim_{n\to\infty} \sum_{i\in\calI} \frac{\mathcal{F}_i^+\bigr(\lambda_n\bigl|\tilde C z_0 \bigr)}{1+\mathcal{F}_i^+\bigr(\lambda_n\bigl|\tilde C  z_0\bigr)}
= \sum_{i\in\calI} \frac{\mathcal{F}_i^+\bigr(\lambda_0\bigl| \tilde C z_0\bigr)}{1+\mathcal{F}_i^+\bigr(\lambda_0\bigl| \tilde C  z_0\bigr)} <1,
\end{eqnarray*}
a contradiction. Consequently, $\lim_{n\to\infty} z^\ast(\lambda_n) = z^\ast(\lambda_0)$.
\par
It remains to prove continuity from the right at $\lambda_0=1$. In this case $\xi_i(1)<1$. Once again $\mathcal{F}_i^+(\lambda_0+\delta|1)<\infty$ for all $\delta>0$ with $\xi_i(\lambda_0+\delta)<1$ according to (\ref{equ:fi+}). Let $(\lambda_n)_{n\in\N}$ be a stricly decreasing sequence of real numbers with limit $1$. We write $z_0=\lim_{n\to\infty} z^\ast(\lambda_n)\leq 1$. Then, for $n$ large enough,
\begin{eqnarray*}
1 & = & \lim_{n\to\infty} \sum_{i\in\calI}
\frac{\mathcal{F}_i^+\bigr(\lambda_n\bigl|z^\ast(\lambda_n)\bigr)}{1+\mathcal{F}_i^+\bigr(\lambda_n\bigl|z^\ast(\lambda_n)\bigr)}\\
&\leq & \lim_{n\to\infty} \sum_{i\in\calI} \frac{\mathcal{F}_i^+\bigr(\lambda_n\bigl|z_0\bigr)}{1+\mathcal{F}_i^+\bigr(\lambda_n\bigl|z_0\bigr)}
= \sum_{i\in\calI} \frac{\mathcal{F}_i^+\bigr(1\bigl|z_0\bigr)}{1+\mathcal{F}_i^+\bigr(1\bigl|z_0\bigr)}.
\end{eqnarray*}
In order to finish the proof we verify that $z^\ast(1)=1$, from which $z_0=z^\ast(1)=1$ follows.
%
Indeed, by Equation (\ref{equ:fi+}) we get
$$
\sum_{i\in\mathcal{I}} \frac{\mathcal{F}_i^+(1|1)}{1+\mathcal{F}_i^+(1|1)} = \sum_{i\in\mathcal{I}}\Bigl(1-G_i\bigl(e_i,e_i|\xi_1(1)\bigr) \bigl(1-\xi_i(1)\bigr) \Bigr).
$$
From \cite[Lemma 5.1]{gilch:07} follows that $1-G_i\bigl(e_i,e_i|\xi_1(1)\bigr) \bigl(1-\xi_i(1)\bigr)$ is just the probability that a single random walk on $\Gamma$ tends to an infinite word of the form $x_1x_2\dots\in\Omega_\infty$ with $x_1\in\Gamma_i^\times$, that is, the above sum equals $1$. 
\end{proof}
The next result completes the proof of Theorem \ref{Th:lambda-function} (2):
\begin{Lemma}
For all $\lambda\in [1,R]$, $\mathrm{HD}(\Lambda)\leq \frac{1}{2}\mathrm{HD}(\Omega)$.
\end{Lemma}
\begin{proof}
Define the function
$$
\mathcal{F}^{(2)}(\lambda|z) := \sum_{x\in\Gamma} F(e,x|\lambda)^2\,z^{l(x)},
$$
whose radius of convergence is denoted by $z^\ast_2$. The Cauchy-Schwarz Inequality gives then
\begin{eqnarray*}
\frac 1{z^{\ast}}&= &\limsup_{m\to\infty}{\Bigl( \sum_{x\in\Gamma: l(x)=m} F(e,x|\lambda)\Bigr)^{1/m}} \\
&\leq &
{\limsup_{m\to\infty}\sqrt{\Bigl( \sum_{x\in\Gamma: l(x)=m} F(e,x|\lambda)^2\Bigr)^{1/m}}}\cdot
\limsup_{m\to\infty}\sqrt{{\Bigl(\sum_{x\in\Gamma: l(x)=m} 1^2\Bigr)^{1/m}}}\\
&=& \sqrt{\frac1{z^\ast_2}} \cdot \sqrt{ \frac1{z^\ast_{\mathcal{S}}}}.
\end{eqnarray*}
To prove the claim of the lemma it suffices (by the formulas given in Theorems \ref{thm:boxdim} and \ref{thm:boxdim-omega}) to show that $z^\ast_2 \geq 1$. First,
\begin{eqnarray*}
\mathcal{F}^{(2)}(\lambda|1) & = & \sum_{x\in\Gamma} F(e,x|\lambda)^2
= \frac{1}{G(e,e|\lambda)^2} \sum_{x\in\Gamma} G(e,x|\lambda)^2 \\
&=& \frac{1}{G(e,e|\lambda)^2} \sum_{x\in\Gamma} \Bigl(\sum_{n\geq 0} p^{(n)}(e,x) \lambda^n \Bigr)^2.
\end{eqnarray*}
For given $x\in\Gamma$, the coefficient of $\lambda^n$ in the inner squared sum can -- by symmetry -- be rewritten as
\begin{equation}\label{equ:square-sum}
 \frac{1}{G(e,e|\lambda)^2} \sum_{m=0}^n p^{(m)}(e,x) p^{(n-m)}(x,e).
\end{equation}
Thus, every path $[x_0=e,x_1,\dots,x_n=e]$ of length $n$ (consisting of $n+1$ vertices) from $e$ to $e$ is counted $n+1$ times, since every $x_i$ can play the role of $x$ in Equation (\ref{equ:square-sum}). That is,
$$
\mathcal{F}^{(2)}(\lambda|z) = \frac{1}{G(e,e|\lambda)^2}\sum_{n\geq 0} p^{(n)}(e,e)\cdot (n+1) \cdot \lambda^n = \frac{\lambda\,G'(e,e|\lambda)}{G(e,e|\lambda)^2}+\frac{1}{G(e,e|\lambda)}.
$$
From this follows  $z^\ast_2 \geq 1$ whenever $\lambda <R$ or
$G'(e,e|R)<\infty$, and thus $\mathrm{HD}(\Lambda)\leq \frac{1}{2}\mathrm{HD}(\Omega)$ for $\lambda< R$. By Lemma \ref{lemma:phi-continuity}, the proposed inequality holds -- due to continuity from the left -- also in the case $\lambda=R$.
\end{proof}


%

In order to prove Theorem \ref{Th:lambda-function} (3) we start with the
following lemma:
\begin{Lemma}
For all $i\in\calI$, $G_i'\bigl(e_i,e_i|\xi_i(R)\bigr)<\infty$.
\end{Lemma}
\begin{proof}
From \cite[Prop. 9.18]{woess} follows $\xi_i(R)\leq R_i$, where $R_i$ is the
radius of convergence of $G_i(e_i,e_i|z)$. If $\xi_i(R)< R_i$ then the claim of
the lemma is obvious. Assume now that $\xi_i(R)= R_i$. Then, by
\cite[Lemma 17.1.(a)]{woess},
$R\,G(e,e|R)=R_i\,G_i(e_i,e_i|R_i)/\alpha_i$. Therefore,
$G_i(e_i,e_i|R_i)<\infty$ since $G(e,e|R)<\infty$ by non-amenability of
$\Gamma$. If $G_i'(e_i,e_i|R_i)=\infty$ would hold, we would get a
contradiction to $\xi_i(R)= R_i$ by \cite[Equ. (9.14), Thm. 9.22, Lemma 17.1.(a)]{woess}.
\end{proof}
Let us remark that $F_i'\bigl(e_i,x\bigl|\xi_i(R)\bigr)=F_i'\bigl(x,e_i\bigl|\xi_i(R)\bigr)<\infty$ for all $x\in\Gamma_\ast^\times$; this
can be easily verified with the help of the inequality
$$
\mu_i^{(n+|x|)}(e_i) \geq \mu_i^{(|x|)}(x) \cdot \Prob\bigl[ Y^{(i)}_n=e_i,\forall m<n:
Y^{(i)}_m\neq e_i\mid Y^{(i)}_0=x\bigr]\quad \textrm{ for all } n\in\N,
$$
where $\bigl(Y_n^{(i)}\bigr)_{n\in\N}$ is a random walk on $\Gamma_i$ governed
by $\mu_i$. We proceed now with expanding the Green function $G(z):=G(e,e|z)$ in a
neighbourhood of $z=R$. By  \cite[Prop. 17.4]{woess} and
\cite[Sec. 3 \& 4]{candellero-gilch}, we have
$$
G(z) = \begin{cases}
G(R) + g_1 \cdot \sqrt{R-z} + \mathbf{o}\bigl(\sqrt{R-z}\bigr), & \textrm{if } G'(R)=\infty,\\
G(R) - G'(R) \cdot (R-z) + \mathbf{o}\bigl(R-z\bigr), & \textrm{if } G'(R)<\infty.
\end{cases}
$$
We write in the following $c:=1/2$, if $G'(R)=\infty$, and $c:=1$ otherwise.
The next aim is to show that the functions $F(e,x|z)$,
$x\in\Gamma\setminus\{e\}$, have the same expansions.
\begin{Lemma}
For all $x\in\Gamma\setminus\{e\}$, there are constants $f_x\neq 0$ such that 
$$
F(e,x|z)=F(e,x|R) + f_x \cdot (R-z)^c + \mathbf{o}\bigl((R-z)^c\bigr).
$$ 
\end{Lemma}
\begin{proof}
We consider the case $c=1$ first. By \cite[Lemma 3.2]{candellero-gilch}, we
have $0<\xi_i'(R)<\infty$, that is, we can write
$$
\xi_i(z)=\xi_i(R) -\xi_i'(R) \cdot (R-z) +\mathbf{o}(R-z).
$$
In the following we write $F_i(e_i,x|z)=\sum_{n\geq 1} f_n(x)z^n$ for
$x\in \Gamma_i^\times$. Therefore,
\begin{equation}\label{equ:Fi-exp}
F(e,x|z)=F_i\bigl(e_i,x|\xi_i(z)\bigr) = \sum_{n\geq 1} f_n(x)\,\bigl(\xi_i(R) -\xi_i'(R) \cdot (R-z) +\mathbf{o}(R-z)\bigr)^n.
\end{equation}
The coeffcient of $(R-z)$ is given by
$$
-\xi_i'(R)\cdot \sum_{n\geq 1}n\cdot f_n(x) \cdot \xi_i(R)^{n-1} =
-\xi_i'(R)F_i'\bigl(e_i,x|\xi_i(R)\bigr)\in (-\infty,0).
$$
Recall that, for $x=x_1\dots x_n\in\Gamma\setminus\{e\}$,
$$
F(e,x_1\dots x_n|z) = \prod_{j=1}^n F_{\tau(x_j)}\bigl(e_{\tau(x_j)},x_j\,\bigl|\, \xi_{\tau(x_j)}(z)\bigr).
$$
Now, plugging the expansion (\ref{equ:Fi-exp}) into the above formula gives us
the coefficient of $(R-z)$:
$$
f_x=\sum_{j=1}^n -\xi_{\tau(x_j)}'(R) F_{\tau(x_j)}'\bigl(e_{\tau(x_j)},x_j\,\bigl|\, \xi_{\tau(x_j)}(R)\bigr)
\prod_{\substack{k=1,\\k\neq j}}^n F_{\tau(x_k)}\bigl(e_{\tau(x_k)},x_k\mid
\xi_{\tau(x_k)}(R)\bigr)\in (-\infty,0).
$$
This yields the claim in the case $c=1$.
\par
We now turn to the case $c=1/2$. By \cite[Equ. (9.20)]{woess}, we have
\begin{equation}\label{equ:Gi-G-equ}
\alpha_i z G(z) = \xi_i(z)\,G_i\bigl(\xi_i(z)\bigr).
\end{equation}
Write $\xi_i(z)=\xi_i(R)+X_i(z)$ with $X_i(R)=0$. Our aim is to show
that $X_i(z)$ is of order $\sqrt{R-z}$, from which we can derive the proposed
expansion of $F(e,x|z)$. We rewrite
(\ref{equ:Gi-G-equ}) as
\begin{eqnarray*}
&&\alpha_i \bigl(R-(R-z)\bigr) \cdot \bigl(G(R) +g_1 \sqrt{R-z}
+\mathbf{o}(\sqrt{R-z})\bigr)\\
&=&
\bigl( \xi_i(R)+X_i(z)\bigr) \cdot \sum_{n\geq 0} \mu_i^{(n)}(e_i) \bigl(\xi_i(R)+X_i(z)\bigr)^n.
\end{eqnarray*}
The constant term on the left hand side of the equation is $\alpha_i\, R\, G(R)$,
which equals the constant term on the right hand side
$\xi_i(R)G_i\bigl(\xi_i(R)\bigr)$ by (\ref{equ:Gi-G-equ}). The coefficient of
$\sqrt{R-z}$ on the left hand side is $\alpha_iR g_1\neq 0$. The coefficient of
$X_1(z)$ on the right hand side is given by
$$
\xi_i(R) G_i'\bigl(e_i,e_i|\xi_i(R)\bigr)+G_i\bigl(e_i,e_i|\xi_i(R)\bigr)\in (0,\infty).
$$
Thus, $X_1(z)\sim \sqrt{R-z}$ as $z\uparrow R$, and therefore
$$
F_i\bigl(e_i,x|\xi_i(z)\bigr) = \sum_{n\geq 1} f_n(x)\,\bigl(\xi_i(R) + \hat\xi_i \cdot \sqrt{R-z} +\mathbf{o}(R-z)\bigr)^n
$$
for some $\hat\xi_i< 0$. The rest follows analogously to the case $c=1$ by replacing
$(R-z)$ with $\sqrt{R-z}$.
\end{proof}
Consider now the following difference for $i\in\calI$:
\begin{eqnarray*}
&& \mathcal{F}_i\bigl(R|z^\ast(R)\bigr) - \mathcal{F}_i\bigl(\lambda
|z^\ast(\lambda)\bigr)\\
&=& \sum_{m\geq 1} z^\ast(R)^m \sum_{\substack{x\in\Gamma_i:\\ l(x)=m}} F(e,x|R) \\
&&\quad -
\sum_{m\geq 1} \Bigl(z^\ast(R) -\bigl(z^\ast(R)-z^\ast(\lambda)\bigr)\Bigr)^m
\sum_{\substack{x\in\Gamma_i:\\ |x|=m}}\Bigl[ F(e,x|R)+ f_x (R-\lambda)^c + \mathbf{o}\bigl((R-\lambda)^c\bigr)\Bigr]\\
&=& \sum_{m\geq 1} z^\ast(R)^m \sum_{\substack{x\in\Gamma_i:\\ l(x)=m}} \bigl( -f_x
(R-\lambda)^c - \mathbf{o}\bigl((R-\lambda)^c\bigr)\bigr) \\
&&\quad + \bigl(z^\ast(R)-z^\ast(\lambda)\bigr) \frac{\partial}{\partial z}
\mathcal{F}_i^+\bigl(\lambda|z^\ast(R)\bigr) + \mathbf{o}\bigl(z^\ast(R)-z^\ast(\lambda)\bigr).
\end{eqnarray*}
Moreover,
\begin{eqnarray}
0 & = &\sum_{i\in\calI}
\frac{\mathcal{F}_i\bigl(R|z^\ast(R)\bigr)}{1+\mathcal{F}_i\bigl(R|z^\ast(R)\bigr)}-
\sum_{i\in\calI}
\frac{\mathcal{F}_i\bigl(\lambda|z^\ast(\lambda)\bigr)}{1+\mathcal{F}_i\bigl(\lambda|z^\ast(\lambda)\bigr)}\nonumber\\
&=& \sum_{i\in\calI} \sum_{n\geq 0}
\Bigl(-\mathcal{F}_i\bigl(\lambda|z^\ast(\lambda)\bigr)\Bigr)^{n+1}-\Bigl(-\mathcal{F}_i\bigl(R|z^\ast(R)\bigr)\Bigr)^{n+1}.\label{equ:decomp1}
\end{eqnarray}
Write
\begin{equation}\label{equ:decomp2}
\Bigl(-\mathcal{F}_i\bigl(\lambda|z^\ast(\lambda)\bigr)\Bigr)^{n+1}-\Bigl(-\mathcal{F}_i\bigl(R|z^\ast(R)\bigr)\Bigr)^{n+1}
=
\Bigl(\mathcal{F}_i\bigl(R|z^\ast(R)\bigr)-\mathcal{F}_i\bigl(\lambda|z^\ast(\lambda)\bigr)\Bigr)
\cdot g_n(\lambda),
\end{equation}
where $g_n(R)\neq 0$ for every $n\in\N$. Plugging the decomposition of
$\mathcal{F}_i\bigl(R|z^\ast(R)\bigr) - \mathcal{F}_i\bigl(\lambda|z^\ast(\lambda)\bigr)$
into (\ref{equ:decomp1}) and comparing all error terms yields in view of (\ref{equ:decomp2}) the following
behaviour:
$$
z^\ast(R)-z^\ast(\lambda) \sim
\begin{cases}
\hat C_1 \cdot (R-\lambda), & \textrm{if } G'(R)<\infty,\\
\hat C_2 \cdot \sqrt{R-\lambda}, & \textrm{if } G'(R)=\infty
\end{cases}
$$
for a suitable constant $\hat C_1$, $\hat C_2$ respectively.
The statement (3) of Theorem \ref{Th:lambda-function} follows now from
$$
\log z^\ast(\lambda)-\log z^\ast(R)= \log \Bigl( 1-\frac{1}{z^\ast(R)}\bigl(z^\ast(R)-z^\ast(\lambda)\bigr)\Bigr)
$$
and by the Taylor expansion of $\log(1-x)$ at $x=0$.

\subsection{Proof of Corollary \ref{Cor:finite-groups}}

In a first step we show the following lemma:

\begin{Lemma}
$$
\overline{\mathrm{BD}^{\mathrm{fin}}}(\Lambda)\leq -\frac{\log\theta}{\log\alpha} \quad
\textrm{ and } \quad
\mathrm{BD}^{\mathrm{fin}}(\Omega)= -\frac{\log\varrho}{\log\alpha}.
$$
\end{Lemma}
\begin{proof}
First, we define the matrices $M_0=\bigl(m_0(i,j)\bigr)_{i,j\in\calI}$ and $D_0=\bigl(d_0(i,j)\bigr)_{i,j\in\calI}$ by
$$
m_0(i,j) :=
\begin{cases}
\mathcal{F}_i^+(\lambda), & \textrm{if } i=j,\\
0, & \textrm{otherwise},
\end{cases}
\quad \textrm{ and } \quad 
d_0(i,j) :=
\begin{cases}
|\Gamma_i|-1, & \textrm{if } i=j,\\
0, & \textrm{otherwise}.
\end{cases}
$$
For $m\in\N$, denote by $\mathcal{H}^{\mathrm{fin}}_m$ the random number of visited words of the form $w_1\dots w_m\in\Gamma$.
Then
\begin{eqnarray*}
\mathbb{E}|\mathcal{H}^{\mathrm{fin}}_m|& \leq & \sum_{x\in\mathcal{H}^{\mathrm{fin}}_m} \mathbb{E}Z_\infty(x) \leq
\sum_{x\in\Gamma:\Vert x\Vert =m} F(e,x|\lambda) =\mathds{1}^T M_0 M^{m-1} \mathds{1} \quad \textrm{ and}\\
\hat S(m) & = & \bigl|\bigl\lbrace x\in\Gamma \,\bigl|\, \Vert x\Vert=m \bigr\rbrace\bigr| = \mathds{1}^T D_0 D^{m-1} \mathds{1}.
\end{eqnarray*}
Let $u\in\mathbb{R}^r$ be an  eigenvector w.r.t. the eigenvalue $\theta$ such that $u\geq
\mathds{1}$. Then:
$$
\mathbb{E}|\mathcal{H}^{\mathrm{fin}}_m|\leq \left(
\begin{array}{c}
\mathcal{F}_1(\lambda)\\
\vdots\\
\mathcal{F}_r(\lambda)
\end{array}
\right)^T
M_1^{m-1} u
\leq 
\left(
\begin{array}{c}
\mathcal{F}_1(\lambda)\\
\vdots\\
\mathcal{F}_r(\lambda)
\end{array}
\right)^T
\theta^{m-1} u.
$$
Thus, $\limsup_{m\to\infty} \bigl(\mathbb{E} \mathcal{H}^{\mathrm{fin}}_m\bigr)^{1/m} \leq \theta$. 
Similarily, one can show that $\lim_{m\to\infty} \hat S(m)^{1/m} =\varrho$ by taking eigenvectors $v_1\geq \mathds{1}$ and $v_2\leq \mathds{1}$. Analogously to the proofs of
Lemma \ref{lemma:upperbound} and Propositions \ref{prop:upperbound}, \ref{prop:boxdim-omega} we obtain the claim. 
\end{proof}
\begin{proof}[Proof of Corollary \ref{Cor:finite-groups}]
First, we remark that we dropped the assumption on symmetry of the $\mu_i$'s in the case of free products of finite groups. This assumption is needed in the general case to ensure $F(e,x|\lambda)<1$. This inequality holds also in the present setting: by  \cite[Equation (9.20)]{woess},
$$
\alpha_i z G(e,e|z) = G_i\bigl(e_i,e_i|\xi_i(z)\bigr) \xi_i(z).
$$
Since $G(e,e|R)<\infty$ and $G_i(e_i,e_i|1)=\infty$, we must have $\xi_i(R)<1$, and consequently,
$$
F(e,x_1\dots x_k|\lambda)  =  \prod_{j=1}^k F_{\tau(x_j)}\bigl(e_{\tau(x_j)},x_j|\xi_{\tau(x_j)}(\lambda)\bigr)
<\prod_{j=1}^k F_{\tau(x_j)}\bigl(e_{\tau(x_j)},x_j|1\bigr)=1.
$$
In order to show that $-\log\theta/\log\alpha$ is a lower bound for $\mathrm{HD}^{\mathrm{fin}}(\Lambda)$, we can follow the reasoning in \cite[Section 6]{lalley2000} or also as in Section \ref{subsec:lower-bound}. Analogously to the proof of Theorem \ref{thm:boxdim-omega} we obtain that $\mathrm{HD}^{\mathrm{fin}}(\Omega)=\mathrm{BD}^{\mathrm{fin}}(\Omega)$. 
\end{proof}

\subsection{Proof of Corollary \ref{Cor:amalgam}}

First, we prove the following lemma:
\begin{Lemma}\label{lemma:amalgam-upper-bound}
$$
\overline{\mathrm{BD}^{(H)}}(\Lambda)\leq -\frac{\log\theta_H}{\log\alpha} \quad \textrm{ and } \quad
\mathrm{BD}^{(H)}(\Omega)= -\frac{\log\varrho_H}{\log\alpha}.
$$ 
\end{Lemma}
\begin{proof}
First, we define the matrices $N_0=\bigl(n_0(i,j)\bigr)_{i,j\in\calI}$ and $D_{0,H}=\bigl(d_{0,H}(i,j)\bigr)_{i,j\in\calI}$ by
$$
n_0(i,j) :=
\begin{cases}
\mathcal{F}_i^{(H)}(\lambda), & \textrm{if } i=j,\\
0, & \textrm{otherwise},
\end{cases}
\quad \textrm{ and } \quad 
d_{0,H}(i,j) :=
\begin{cases}
[\Gamma_i:H_i]-1, & \textrm{if } i=j,\\
0, & \textrm{otherwise}.
\end{cases}
$$
For $m\in\N$, we denote by $\mathcal{H}_m^{(H)}$ the set of words of the form $g_1\dots g_mh\in\Gamma$ in the sense of (\ref{def:amalgam-words}). Since every path from $e$ to  $g_1\dots g_mh\in\Gamma$ has to pass through points $g_1\dots g_jh_j\in\Gamma$, where $h_j\in H$ with $h_m=h$, we have
$$
\sum_{g_1\dots g_mh \in\Gamma} F_H(g_1\dots g_m h|z) =
\sum_{g_1\dots g_mh \in\Gamma} \sum_{h_1,\dots,h_{m-1}\in H} \prod_{i=1}^m F_H(g_ih_i|z)=\mathds{1}^T N_0 N^{m-1} \mathds{1}.
$$
Choose now an eigenvector $v=(v_1,\dots,v_r)^T\geq \mathds{1}$ w.r.t. the eigenvalue $\theta_H$ of $N$. Then
$$
\mathbb{E}|\mathcal{H}^{(H)}_m| \leq \mathds{1}^T N_0 N^{m-1} \mathds{1} \leq \mathds{1}^T N_0 N^{m-1} v = \theta_H^{m-1} \cdot \Bigl(\sum_{i\in\calI} v_i \mathcal{F}_i^{(H)}(\lambda) \Bigr),
$$
and therefore, $\limsup_{m\to\infty} \mathbb{E}|\mathcal{H}^{(H)}_m|^{1/m} \leq \theta_H$. 
 Furthermore, we remark that  $\hat S_H(m)=\bigl|\{x_1\dots x_m\mid x_i\in\bigcup_{j\in\calI} \mathcal{R}_j\setminus\{e_j\}, x_i\in\mathcal{R}_j\Rightarrow x_{i+1}\notin \mathcal{R}_j\}\bigr|$ can be written as
$$
\hat S_H(m) = \mathds{1}^T D_{0,H} D_H^{m-1} \mathds{1}.
$$ 
Taking eigenvectors $v_1\geq \mathds{1}$ and $v_2\leq \mathds{1}$
w.r.t. $\varrho_H$ leads to $\lim_{m\to\infty}|\hat S_H(m)|^{1/m}=\varrho_H$. The same reasoning as used in the proofs of
Lemma \ref{lemma:upperbound} and Propositions \ref{prop:upperbound},\ref{prop:boxdim-omega} yields the proposed claim. 
\end{proof}
\begin{proof}[Proof of Corollary \ref{Cor:amalgam}]
It is sufficient to show that $-\log\theta_H/\log\alpha$ is also a lower bound for 
 $\mathrm{HD}^{(H)}(\Lambda)$. First, we remark that for $m\in\N$
\begin{eqnarray*}
\sum_{g_1\dots g_mh\in\Gamma: g_1\notin \mathcal{R}_1} \mathbb{E}Z_\infty(g_1\dots g_mh) & = &
\sum_{g_1\dots g_mh\in\Gamma: g_1\notin \mathcal{R}_1} F(e,g_1\dots g_mh|\lambda) \\
&=& \sum_{g_1\dots g_mh\in\Gamma: g_1\notin \mathcal{R}_1} \sum_{h_0\in H }F_H(g_1\dots g_mh_0|\lambda)\,F(e,h_0^{-1}h|\lambda).
\end{eqnarray*}
Since $|H|<\infty$, there are constants $d,D>0$ such that $d\leq F(e,h|\lambda)\leq D$ for all $h\in H$. We write $\mathds{1}_0:=(0,1,\dots,1)^T\in\mathbb{R}^{r}$ and get:
\begin{eqnarray*}
\Bigl( \sum_{g_1\dots g_mh\in\Gamma: g_1\notin \mathcal{R}_1} \mathbb{E}Z_\infty(g_1\dots g_mh) \Bigr)^{1/m}
& \leq &  \Bigl(D\cdot \mathds{1}_0^T N_0 N^{m-1}\mathds{1} \Bigr)^{1/m}
\xrightarrow{m\to\infty} \theta_H\quad \textrm{ and }\\
\Bigl( \sum_{g_1\dots g_mh\in\Gamma: g_1\notin \mathcal{R}_1}
\mathbb{E}Z_\infty(g_1\dots g_mh) \Bigr)^{1/m}
&\geq &  \Bigl(d\cdot \mathds{1}_0^T N_0 N^{m-1}\mathds{1} \Bigr)^{1/m}
\xrightarrow{m\to\infty} \theta_H.
\end{eqnarray*}
This can be easily verified by substituting $\mathds{1}$ by an eigenvector $v_1\geq \mathds{1}$ of $\theta_H$, by an eigenvector $v_2\leq \mathds{1}$ of $\theta_H$ respectively.
With the help of this convergence behaviour and the last lemma, we can prove once again analogously to the reasoning in \cite[Section 6]{lalley2000} or Section \ref{subsec:lower-bound} that the upper bounds in Lemma \ref{lemma:amalgam-upper-bound} equal the Hausdorff and the Box-Counting dimensions. Analogously to the proof of Theorem \ref{thm:boxdim-omega} we obtain that $\mathrm{HD}^{(H)}(\Omega)=\mathrm{BD}^{(H)}(\Omega)$.
\end{proof}

\section*{Acknowledgements}
The authors are grateful to the referee whose comments led to an essential improvement of the paper.

The research of Elisabetta Candellero was financially supported by the Austrian
Academy of Science (\"OAW), the research of Lorenz Gilch by the German Research
Foundation (DFG) grant GI 746/1-1, while Sebastian M\"uller was supported  by
PIEF-GA-2009-235688. Part of the work was done during a visit of S. M\"uller
at Graz University of Technology and a visit of E. Candellero at Geneva University, where both visits were supported by the ESF Grant ``Random Geometry of Large Interacting Systems and Statistical Physics''. The authors are also grateful to Steven Lalley for discussions during his visit at Graz University of Technology.

\bibliographystyle{abbrv}
\bibliography{literatur}

\end{document}